\documentclass[sn-basic,smallextended]{sn-jnl}

\usepackage[binary-units=true]{siunitx}
\usepackage{amsmath,amssymb}
\usepackage{amsthm}
\usepackage{mathtools}
\usepackage{array} 
\usepackage{varwidth} 
\usepackage{array}
\usepackage{multirow}

\usepackage{nomencl}
\usepackage{xpatch}
\xpatchcmd{\thenomenclature}{\section*{\nomname}}{}{\typeout{Success}}{\typeout{Failure}}
\xpatchcmd{\thenomenclature}{\section*{Nomenclature}}{}{\typeout{Success}}{\typeout{Failure}}

\usepackage{color, colortbl}
\definecolor{mygray}{gray}{0.6}

\allowdisplaybreaks

\makenomenclature

\DeclareMathOperator{\diag}{diag}

\DeclareMathOperator*{\mini}{min.}
\DeclareMathOperator*{\lexmin}{lexmin.}

\definecolor{light-gray}{gray}{0.95}
\definecolor{dark-gray}{gray}{0.5}
\definecolor{mygray}{gray}{0.75}
\newcommand{\BIN}{\begin{bmatrix}}
\newcommand{\BOUT}{\end{bmatrix}}

\newcommand{\tA}{\tilde{A}}
\newcommand{\utA}{\underline{\tilde{A}}}

\newcommand{\mA}{{\mathcal{A}}}
\newcommand{\umI}{\underline{\mathcal{I}}}

\newcommand{\mI}{{\mathcal{I}}}
\newcommand{\bE}{{\mathbb{E}}}
\newcommand{\bI}{{\mathbb{I}}}

\definecolor{orange}{rgb}{0.99,0.69,0.07}
\definecolor{lightgray}{gray}{0.85}
\definecolor{light-gray}{gray}{0.95}
\definecolor{dark-gray}{gray}{0.5}

 \makeatletter
 \newcommand\fs@spaceruled{\def\@fs@cfont{\bfseries}\let\@fs@capt\floatc@ruled
   \def\@fs@pre{\vspace{5pt}\hrule height.8pt depth0pt \kern2pt}%
   \def\@fs@post{\kern2pt\hrule\relax}%
   \def\@fs@mid{\kern2pt\hrule\kern2pt}%
   \let\@fs@iftopcapt\iftrue}
 \makeatother

\newtheorem{theorem}{Theorem}

\title[$\mathcal{N}$\hspace{-2pt}IPM-HLSP]{$\mathcal{N}$\hspace{-0.5pt}IPM-HLSP: An Efficient Interior-Point Method for Hierarchical Least-Squares Programs}
\author*[1]{\fnm{Kai} \sur{Pfeiffer}}\email{kaipfeifferrobotics@gmail.com}
\author[2]{\fnm{Adrien} \sur{Escande}}
\author[3]{\fnm{Ludovic} \sur{Righetti}}
\affil*[1]{\orgdiv{School of Mechanical and Aerospace Engineering}, \orgname{Nanyang Technological University}, \country{Singapore}}
\affil[2]{\orgdiv{Inria Center of Grenoble Alpes University}, \orgname{Montbonnot}, \country{France}}
\affil[3]{\orgdiv{Tandon School of Engineering}, \orgname{New York University}, \country{New York, USA}}

\begin{document}

	\abstract{
		Hierarchical least-squares programs with linear constraints (HLSP) are a type of optimization problem very common in robotics. Each priority level contains an objective in least-squares form which is subject to the linear constraints of the higher priority levels. Active-set methods are a popular choice for solving them. However, they can perform poorly in terms of computational time if there are large changes of the active set. We therefore propose a computationally efficient primal-dual interior-point method (IPM) for dense HLSP's which is able to maintain constant numbers of solver iterations in these situations. We base our IPM on the computationally efficient nullspace method as it requires only a single matrix factorization per solver iteration instead of two as it is the case for other IPM formulations. We show that the resulting normal equations can be expressed in least-squares form. This avoids the formation of the quadratic Lagrangian Hessian and can possibly maintain high levels of sparsity. Our solver reliably solves ill-posed  instantaneous hierarchical robot control problems without exhibiting the large variations in computation time seen in active-set methods.
	}

	\maketitle 
	
	\printkeywords

	\section{Introduction} 
	\subsection{Context and contribution}
	
	Hierarchical or lexicographic multi-objective optimization (LMOO) is the prioritized separation of constraints and objectives over any number of priority levels $p$. 
	In its most generic form, LMOO can be written as ($\lexmin$: lexicographically minimize)
	\begin{align}
		\lexmin_u\qquad &f_{\mathbb{E}_1}(u), \dots , f_{\mathbb{E}_p}(u)
		\label{eq:HNLP}\tag{LMOO}\\
		\text{s.t}\qquad& f_{\mathbb{I}_{\cup p}}(u) \geq 0\nonumber
	\end{align}
	$u\in\mathbb{R}^n$ is a variable vector. 
	Each (possibly) non-linear \textit{objective}  $f_{\mathbb{E}_l}\in\mathbb{R}^{m_{\mathbb{E}_l}}$, which represent the set of equality constraints $\mathbb{E}_l$ of level $l$ ($\lvert\mathbb{E}_l\rvert = m_{\mathbb{E}_l}$), needs to be minimized. This comes to the extent of not influencing the optimality of the infinitely more important \textit{constraints} $f_{\mathbb{E}_{\cup l-1}}\in\mathbb{R}^{m_{\mathbb{E}_{\cup l-1}}}$ of levels $1$ to $l-1$. $\mathbb{E}_{\cup l-1}$ represents the set union $\mathbb{E}_{\cup l-1} \coloneqq \bigcup_{i=1}^{l-1}\mathbb{E}_i = \mathbb{E}_1 \cup \cdots \cup \mathbb{E}_{l-1}$ of the sets of equality constraints of levels 1 to $l-1$. At the same time, inequality constraints $f_{\mathbb{I}_{\cup p}}\in\mathbb{R}^{m_{\mathbb{I}_{\cup p}}}$  shall not be violated.

	To date, no method has been proposed to efficently solve \ref{eq:HNLP} as is. Early beginnings of LMOO have been reported for example in~\cite{Sherali1983} in the context of lexicographic linear programming. The authors propose a set of appropriate weights to solve the prioritized multi-objectives as an equivalent weighted one. A solution in a variable reducing fashion has been proposed in~\cite{Holder2006}. In this case the problem dimensions decrease as the solver progresses through the lexicographical program. This is one of the fundamental ideas in efficient hierarchical programming. Extensions to non-linear objective functions based on optimal Pareto front and for lexicographic programs and non-smooth functions can be found in~\cite{Evtushenko2014},~\cite{Lai2021} and~\cite{abu2023}, respectively. A solution can also be obtained by evolutionary algorithms with great efficiency in identifying trade-off solutions. An evaluation on the scheduling problem for satellite communication is given in~\cite{Petelin2021}. 
	
	A specific form of LMOO is hierarchical linear least-squares programming (HLSP), which is the main focus of this work. Here all constraints are linear and the objectives are in least-squares form. HLSP's have seen a sharp rise in popularity in the robotics community over the recent years. Different works tackled the incorporation of infeasible inequality constraints on any priority level~\citep{Kanoun2009}, or proposed very efficient solvers for these types of problems~\citep{escande2014}. The aforementioned HLSP solver~\citep{escande2014} enables high frequency non-linear whole-body robot control when used within a sequential hierarchical least-squares programming (S-HLSP) solver to solve non-linear hierarchical least-squares programs (\ref{eq:NL-HLSP}). S-HLSP with local convergence properties for example based on a real-time capable trust-region adaptation method~\citep{pfeiffer2018} has been applied in many robot control works (\cite{Herzog2014} and references therein). In contrast, a globally converging method for non-linear least-squares programming based on branch-and-bound has been proposed in~\cite{Sahinidis2012}. However, this method is applicable only to unconstrained problems.
	
	In this work we propose an efficient, dense HLSP solver based on the interior point method (IPM). 
	This has numerical advantages with respect to active-set method (ASM) based solvers like~\cite{escande2014}. The \textit{active-set method} separates inequality constraints into \textit{active} constraints, which are violated, and \textit{inactive} constraints, which are satisfied. Violated constraints are commonly referred to as \textit{infeasible}. In ill-posed robot control scenarios our IPM based solver maintains constant number of iterations and computation times as opposed to the ASM. This is favorable in real-time robot control where only limited computational resources are available and a control loop has a fixed time budget at each iteration. 
	
	The solver's C++ code based on the Eigen library~\citep{eigenweb} is available at \url{https://www.github.com/pfeiffer-kai/NIPM-HLSP}.
	The nomenclature and variable naming used hereafter is listed below.

	\printnomenclature[2cm]
	
	\section*{Nomenclature}
	\begin{itemize}[align=parleft]
		\item[$l$] Current priority level
		\item[$l^*$] Virtual priority level
		\item[{$p$}] Overall number of priority levels, excluding the trust region constraint on $l=0$
		\item [{$n$}] Number of variables
		\item [{$r$}] Rank of matrix 
		\item [{$n_r$}] Number of remaining variables after nullspace projections 
		\item [{$m$}] Number of constraints 
		\item [{$x\in\mathbb{R}^{n}$}] Primal of HLSP 
		\item [{$\Delta x\in\mathbb{R}^{n}$}] Primal Newton step of HLSP 
		\item [{$\Delta z$}] Primal nullspace step 
		\item [{$f(u)\in\mathbb{R}^{m}$}] Non-linear constraint function of variable vector $u\in\mathbb{R}^n$ 
		\item [{$\mathbb{E}_l$}] Set of $m_{\mathbb{E}}$ equality constraints (eq.)   of level $l$
		\item [{$\mathbb{I}_l$}] Set of $m_{\mathbb{I}}$ inequality constraints (eq.)  of level $l$
		\item [{$\mathcal{I}_l$}] Set of $m_{\mathcal{I}}$ inactive inequality constraints (ineq.)  of level $l$ 
		\item [{$\mathcal{A}_l$}] Set of $m_{\mathcal{A}}$ active equality and inequality constraints of level $l$
		\item [{$\mathbb{E}_{\cup l}$ (or ${\mathcal{E}}_{\cup l}$)}] Set union $\mathbb{E}_{\cup l}\coloneqq\bigcup_{i=1}^l \mathbb{E}_i= \mathbb{E}_1 \cup \cdots \cup \mathbb{E}_l$ with $m_{{\mathbb{E}}_{\cup l}}$ constraints
		\item [{$A_{\mathbb{E}}\in\mathbb{R}^{m_{\mathbb{E}}\times n}$}] Matrix  representing a set $\mathbb{E}$  of $m_{\mathbb{E}}$ linear constraints 
		\item [{$b_{\mathbb{E}}\in\mathbb{R}^{m_{\mathbb{E}}}$}] Vector representing a set $\mathbb{E}$ of $m_{\mathbb{E}}$ linear constraints 
		\item [{$\mathcal{N}(A_{\mathcal{A}_{l}})$}] Operator to compute the nullspace basis $Z_{\mA_{l}}$ and the rank $r$ of a matrix $A_{\mathcal{A}_{l}}$
		\item [{$Z_{\mA_{l}}\in\mathbb{R}^{n\times n_r}$}] Nullspace basis of matrix $A_{\mathcal{A}_{l}}\in\mathbb{R}^{m_{\mathcal{A}_{l}}\times n}$ with rank $r$, $n_r = n-r$ and $A_{\mathcal{A}_{l}}Z_{\mA_{l}}=0$ 
		\item [{$N_{\mathcal{A}_{\cup l}}\in\mathbb{R}^{n\times n_r}$}] Accumulated nullspace basis $N_{\mathcal{A}_{\cup l}} = Z_{\mA_{1}} \dots Z_{\mA_{l}}$ 
		\item [{$\tilde{M}\in\mathbb{R}^{m\times n_r}$}] Matrix $\tilde{M} = MN$ projected into the nullspace basis $N\in\mathbb{R}^{n\times n_r}$ of a matrix $A\in\mathbb{R}^{m\times n}$ of rank $r$; $n_r=n-r$  (variable elimination)
		\item [{$v\in\mathbb{R}^{m}$}] Slack variable
		\item [{$V\in\mathbb{R}^{{m,m}}$}] Diagonal matrix equivalent $V=\diag(v)$ of vector $v\in\mathbb{R}^m$  
		\item [{$v^*\in\mathbb{R}^{m}$}] Optimal slack variable 
		\item [{$\lambda\in\mathbb{R}^{{m}}$}] Lagrange multiplier  
		\item [{$\mathcal{L}$}] Lagrangian 
		\item [{$K$}] Gradient of Lagrangian $K \coloneqq \nabla \mathcal{L}$
		\item [{$\iota$}] Newton iteration or active-set iteration of HLSP solver 
		\item [{$\mu$}] Duality measure 
		\item [{$\sigma$}] Centering parameter 
		\item [{$\rho$}] Trust region radius 
		\item [{$\xi$}] Activation threshold of inequality constraints  
		\item [{$\epsilon$}] Convergence threshold for Newton's method 
		\item [{$e\in\mathbb{R}^m$}] Vector of ones  
		\item [{$a\odot b$}] Element-wise multiplication between two vectors $a$ and $b$
	\end{itemize}
	
	\nomenclature[]{$l$}{Current priority level}
	\nomenclature[]{$p$}{Overall number of priority levels, excluding the trust region constraint on $l=0$}
	\nomenclature[]{$n$}{Number of variables}
	\nomenclature[]{$r$}{Rank of matrix} 
	\nomenclature[]{$n_r$}{Number of remaining variables after nullspace projections} 
	\nomenclature[]{$m$}{Number of constraints} 
	\nomenclature[]{$m_l$}{Sum of number of constraints of levels up to $l$} 
	\nomenclature[]{$x\in\mathbb{R}^{n}$}{Primal of HLSP} 
	\nomenclature[]{$\Delta x\in\mathbb{R}^{n}$}{Primal Newton step of HLSP} 
	\nomenclature[]{$\Delta z$}{Primal nullspace step} 
	\nomenclature[]{$f(u)\in\mathbb{R}^{m}$}{Non-linear task function of variable vector $u\in\mathbb{R}^n$} 
	\nomenclature[]{$\mathbb{E}$}{Set of equality constraints} 
	\nomenclature[]{$\mathcal{I}$}{Set of inactive inequality constraints} 
	\nomenclature[]{$\mathcal{A}$}{Set of active equality and inequality constraints} 
	\nomenclature[]{$\underline{\mathcal{A}}_l$}{Set union $\mathcal{A}_1 \cup \cdots \cup \mathcal{A}_l$} 
	\nomenclature[]{$A\in\mathbb{R}^{m\times n}$}{Constraint matrix} 
	\nomenclature[]{$A_{\mathcal{A}}$}{Matrix $A$ (or vector)  representing the constraint set $\mathcal{A}$}  
	\nomenclature[]{$b\in\mathbb{R}^{m}$}{Constraint right hand side} 
	\nomenclature[]{$\underline{A}_i\in\mathbb{R}^{\sum_i^l m_i\times n}$}{Stacked matrix $ \begin{bmatrix} {A}_1^T & \cdots & {A}_l^T \end{bmatrix}^T $} 
	\nomenclature[]{$\underline{b}_i\in\mathbb{R}^{\sum_i^l m_i}$}{Stacked vector $ \begin{bmatrix} {b}_1^T & \cdots & {b}_l^T \end{bmatrix}^T $} 
	\nomenclature[]{$\mathcal{N}(A)$}{Operator to compute the nullspace basis and the rank of a matrix $A$}
	\nomenclature[]{$Z\in\mathbb{R}^{n\times n_r}$}{Nullspace basis of matrix $A\in\mathbb{R}^{m\times n}$ with rank $r$ and $AZ=0$} 
	\nomenclature[]{$N_l\in\mathbb{R}^{n\times n_r}$}{Accumulated nullspace basis $N_l = Z_1 \dots Z_l$} 
	\nomenclature[]{$\tilde{M}\in\mathbb{R}^{m\times n_r}$}{Matrix $\tilde{M} = MN$ projected into the nullspace basis $N\in\mathbb{R}^{n\times n_r}$ of a matrix $A\in\mathbb{R}^{m\times n}$ of rank $r$; $n_r=n-r$} 
	\nomenclature[]{$v_l\in\mathbb{R}^{m}$}{Slack variable of level $l$; we have $v_l = \Lambda_{l,l}$} 
	\nomenclature[]{$v^*\in\mathbb{R}^{m}$}{Optimal slack variable} 
	\nomenclature[]{$\lambda\in\mathbb{R}^{{m}}$}{Lagrange multiplier}  
	\nomenclature[]{$V\in\mathbb{R}^{{m,m}}$}{Diagonal matrix equivalent $V=\diag(v)$ of vector $v\in\mathbb{R}^m$}  
	\nomenclature[]{$\mathcal{L}$}{Lagrangian} 
	\nomenclature[]{$K$}{Gradient of Lagrangian} 
	\nomenclature[]{$\iota$}{Newton iteration or active-set iteration of HLSP solver} 
	\nomenclature[]{$\mu$}{Duality measure} 
	\nomenclature[]{$\sigma$}{Centering parameter} 
	\nomenclature[]{$\rho$}{Trust region radius} 
	\nomenclature[]{$\xi$}{Activation threshold of inequality constraints}  
	\nomenclature[]{$\epsilon$}{Convergence threshold for Newton's method} 
	\nomenclature[]{$e\in\mathbb{R}^m$}{Vector of ones}  
	\nomenclature[]{$\odot$}{Element-wise multiplication between two vectors}  
	
	\subsection{Hierarchical least-squares programming}
	
	A HLSP is a sub-form of~\ref{eq:HNLP} with linear constraints and least-squares objectives as follows
	\begin{align}
		\lexmin_{x,v_{\bE_{\cup p}},v_{\bI_{\cup p}}}\qquad &\Vert v_{\mathbb{E}_1}\Vert^2 + \Vert v_{\mathbb{I}_1}\Vert^2, \dots , \Vert v_{\mathbb{E}_p}\Vert^2 + \Vert v_{\mathbb{I}_p}\Vert^2 \label{eq:hlsplexmin}\\
		\text{s.t}\qquad& A_{\bE_{\cup p}}x -b_{\bE_{\cup p}} = v_{\bE_{\cup p}}\nonumber\\
		& A_{\bI_{\cup p}}x -b_{\bI_{\cup p}} \geq v_{\bI_{\cup p}}\nonumber
	\end{align}
	The variable vector $x\in\mathbb{R}^n$ consists of $n$ variables. The linear constraint sets $\mathbb{E}_l$ and $\mathbb{I}_l$ of each level $l=1,\dots,p$ are represented by the constraint matrices and vectors $A_{\mathbb{E}_l}\in\mathbb{R}^{m_{\mathbb{E}_l}}$ and $b_{\mathbb{E}_l}\in\mathbb{R}^{m_{\mathbb{E}_l}}$, respectively.
	The slack variables $v_{\bE_l}\in\mathbb{R}^{m_{\mathbb{E}_l}}$ and $v_{\bI_l}\in\mathbb{R}^{m_{\mathbb{I}_l}}$ relax the equality and inequality constraints ${\bE_l}$ and ${\bI_l}$ in case of infeasibility. This relaxation allows handling of infeasible inequality constraints on any priority level~\citep{Kanoun2011}. 
	In contrary, LMOO as treated in~\cite{Lai2021} only allows the incorporation of inequality constraints on the problem variables. On each level $l=1,\dots,p$ the optimal slacks $v_{\mathbb{E}_l}^*$ and $v_{\mathbb{I}_l}^*$ need to be identified. At the same time the optimal slacks $v_{\mathbb{E}_{\cup l-1}}^*$ and $v_{\mathbb{I}_{\cup l-1}}^*$ of the previous levels $1$ to $l-1$ must not be violated. 
	
	In essence, a HLSP with $p$ levels is an optimization problem composed of $p$ consecutive least-squares programs (LSP). Each level $l$ consists of a least squares objective which is subject to the linear constraints associated with the higher priority levels $1$ to $l-1$ (except for the first level 1 which does not carry any constraints from previous levels).
	The corresponding LSP's of levels $l=1,\dots, p$ can be written as
	\begin{equation}
		\begin{aligned}
			\mini_{x,v_{\mathbb{E}_l},v_{\mathbb{I}_l}}& \qquad \frac{1}{2}\Vert v_{\mathbb{E}_l} \Vert^2 + \frac{1}{2}\Vert v_{\mathbb{I}_l} \Vert^2\qquad l=1,\dots,p\\
			\text{s.t.}
			& \qquad A_{\mathbb{E}_l}x - b_{\mathbb{E}_l} = v_{\mathbb{E}_l}\\
			& \qquad A_{\mathbb{I}_l}x - b_{\mathbb{I}_l} \geq v_{\mathbb{I}_l}\\
			& \qquad A_{\bE_{\cup l-1}}x - b_{\bE_{\cup l-1}} = v_{\bE_{\cup l-1}}^*\\
			& \qquad A_{\bI_{\cup l-1}}x - b_{\bI_{\cup l-1}} \geq 0
		\end{aligned}
		\label{eq:hlspLvll}
	\end{equation}
	$v_{\bE_{\cup l-1}}^*$ indicates the optimal slacks that were identified for the equality constraints of the previous levels $1$ to $l-1$. Inequality constraints are assumed to be feasible on each priority level. In Sec.~\ref{sec:ipmHLSP} we introduce a relaxation of this assumption based on the active-set method such that the original lexicographical problem~\eqref{eq:hlsplexmin} is represented exactly. 
	
	The authors in~\cite{Kanoun2009} consecutively solve each level $l$
	from the first to the last level $p$ of~\eqref{eq:hlspLvll} in a cascade-like manner. The disadvantage of this approach is that each optimization problem grows in the number of constraints since the ones from the previous levels need to be carried over. In~\cite{DeLasa2009} this is avoided by solving each level in the nullspace of the active constraints of the previous levels. This way the problem dimensions are progressively reduced. This nullspace projection and consequent variable reduction is commonly referred to as the \textit{nullspace method}~\citep{Nocedal2006}. Additionally, only feasible inequality constraints are considered on the first level.  The approach considered in~\cite{escande2014} is based on this same principle but enables the incorporation of inequality constraints on any priority level. Unlike the previous cascade-like approaches~\citep{Kanoun2009,DeLasa2009}, their active-set search solves the whole hierarchy at once which adds further efficiency by unifying the active-set search of all the levels into one. 
	
	The dedicated HLSP solver in~\cite{escande2014} is based on the ASM. At each solver iteration an equality only problem consisting of the active constraints $\mA_{\cup p}$ of the priority levels $1$ to $p$ is solved. $\mA_{\cup p}$ assembles equality constraints $\bE_{\cup p}$ and active inequality constraints of $\bI_{\cup p}$. Based on the violation or feasibility of inequality constraints  $\bI_{\cup p}$, the active set $\mA_{\cup p}$ is composed accordingly. The ASM converges if no constraint needs to be added to or removed from the active set. In ASM's based on the nullspace method, all constraint matrices of levels $l$ to $p$ are projected into the nullspace of the active constraints ${\mA_{\cup l-1}}$ of levels 1 to $l-1$~\citep{coleman1984,benzi2005}. With the right choice of nullspace basis this leads to a possibly significant reduction of variables especially on lower priority levels in case of a large number of linearly independent active constraints. This makes the solvers very efficient as the matrix factorization for the linear system solution is cubically dependent of the number of variables. While the work in~\cite{escande2014} uses orthogonal nullspace bases, the authors in~\cite{dimitrov:2015} implement non-orthogonal bases with further computational advantage due to their block structure. The disadvantage of non-orthogonal bases is that the primal $x$ is not of minimum norm but this characteristic can be easily enforced by regularization.
	
	ASM based solvers are most efficient on problems with limited changes of the active-set between HLSP problem instances. 
	This is usually the case for parametric problems like robot control scenarios where the problem variables evolve slowly. In this case the ASM can be warm started by setting the active set of the current problem with the active set from the previous one~\citep{gill1986lssol}.  However, in robot real-time control there are cases where large shifts of the active set occur, for example in instances close to contact shifting or oscillations in case of numerical instabilities due to ill-posed constraint matrices~\citep{pfeiffer2018,pfeiffer2023}.  In these cases the number of iterations of the ASM may be exponential in the number of constraints as all possible active-set combinations need to be explored~\citep{Rao1998}. Numerical issues like cycling (repeated activation and deactivation of the same constraint) can further increase the number of active set iterations such that it becomes impractical for real-time robot control. While methods mitigating such effects like decomposition updates~\citep{hammarling2008} and cycling handling~\citep{gill1989} exist, they might not be enough to make up for the possibly large number of active set iterations until convergence~\citep{pfeiffer2023}.
	
	In contrast, the IPM has been shown to robustly converge in a deterministic number of iterations independent of problem conditioning~\citep{Nesterov1994,bartlett2000}. The IPM has been first developed for linear programming~\citep{Karmarkar1984} but has seen extensions to  quadratic programs (which LSP's are a sub-form of,~\cite{vanderbei1999}) or non-linear programming~\citep{ipopt}.
	The ASM only considers inequality constraints which are deemed active in the current guess of the active-set. In contrast, the IPM considers all constraints including all inactive inequality constraints. In primal-dual formulations of the IPM, the primal and dual variables thereby move within the interior of the feasible region. This is first achieved by maintaining feasibility conditions by line search or more complex methods based on non-linear arc-searches~\citep{Yang2022}.
	Furthermore, violations are penalized for example by a log-barrier function. A function with an upper bound on solver iterations has been proposed for linear programs  in~\cite{Fathi2020}. A summary of different functions is given in~\cite{Li2021}. 
	A single iteration of the IPM is relatively expensive but convergence is achieved robustly. The authors in~\cite{kuindersma:icra:2014} use this fact to switch from the ASM to the reliable IPM in case of ASM failures.

	\subsection{Our contribution}
	
	While the IPM has been developed for LSP~\citep{vanderbei1999} and could be used to solve~\eqref{eq:hlspLvll} directly by solving each level in sequence~\citep{Kanoun2009}, there exist no dedicated IPM based solvers for HLSP. In this work we provide the theoretical foundations for an efficient IPM for HLSP. 
	The IPM based HLSP solver provides predictability in terms of computation times due to its constant number of Newton iterations even in case of ill-posed constraints matrices. This can be important for example when critical safety constraints need to be dealt with but the ASM fails to find an optimal point in a reasonable number of active set iterations.
	
	Our contributions are threefold: 
	\begin{itemize}
		\item We formulate the IPM for HLSP. It reliably resolves ill-posed HLSP's without significant fluctuations in solver iterations or computation time. We show that this enables a humanoid robot to handle large changes of the active set after a strong push.
		\item We suggest to apply the nullspace method~\citep{Nocedal2006} to the IPM for HLSP instead of the commonly used Schur complement~\citep{wangboyd2010,zeilinger2012}. This transfers the same variable reducing quality to the IPM for HLSP as seen for the ASM for HLSP. While this has been previously done for QP's in Model predictive control (MPC)~\citep{HPIPM}, we give a more detailed explanation for example on our choice of nullspace basis. With this formulation we reduce the necessary number of decompositions of the KKT system per Newton iteration from two to one. We show that this is computationally efficient for almost all problem constellations. Furthermore, the Lagrange multipliers associated with active constraints do not need to be evaluated by virtue of an adapted IPM convergence test.
		\item We show that the IPM can be expressed in least squares form. This way the formation of expensive matrix products is avoided which is of advantage for problems with high number of variables but low number of constraints. This also potentially preserves a large degree of sparsity of the constraint matrices which we aim to exploit in a future sparse version of the proposed solver.
	\end{itemize}  
	
	\subsection{Overview}
	
	This article is structured as follows: In Sec.~\ref{sec:ipmHLSP} we first introduce the notion of active sets into the optimization problem~\eqref{eq:hlspLvll}. We then continue to outline the formulation of the IPM for HLSP. Section~\ref{sec:compHLSP} oversees the efficient computation of the IPM solver iterations. First, the overall algorithm is outlined in Sec.~\ref{sec:algo}. We apply the nullspace method which is a common tool in hierarchical programming, see Sec.~\ref{sec:nsproj}. This requires a concept we refer to as `virtual priority level' and which is detailed in Sec.~\ref{sec:virtualPriorLvl}. With a special convergence test the calculation of the dual associated with the active constraints  can be avoided (Sec.~\ref{sec:ipmnodual}). In Sec.~\ref{sec:lqipm} we show that the solver iterations can be expressed in efficient least-squares form.
	Section~\ref{sec:predcoralgo} oversees the development of Mehrotra's predictor-corrector-algorithm for HLSP's. 
	Finally, we give a computational comparison between the different solver formulations (Sec.~\ref{sec:qpopcomp}).
	The proposed algorithms are evaluated in Sec.~\ref{sec:eval}. We conclude the article with some remarks and considerations for future work (see Sec.~\ref{sec:conclusion}).

	\section{Formulating the interior point method for HLSP}
	\label{sec:ipmHLSP}

	We first separately consider the top three lines of the optimization problem~\eqref{eq:hlspLvll}.
	This problem corresponds to finding a \textit{feasible point} ($v_{\mathbb{E}_l}^*=0$, $v_{\mathbb{I}_l}^*\geq 0$) or \textit{optimal infeasible point} ($v_{\mathbb{E}_l}^*\neq0$, $v_{\mathbb{I}_l}^* < 0$ with minimal $ \Vert v_{\mathbb{E}_l}\Vert^2$, $\ \Vert v_{\mathbb{I}_l}\Vert^2$) of the constraints. 
	In the context of linear programming the self-dual embedding model~\citep{ye1994} has been proposed. Similarly to~\cite{wangboyd2010} for quadratic programming, infeasible initial points with respect to the equality constraints are handled but inequality constraints are assumed to be feasible with $ A_{\mathbb{I}_l}x - b_{\mathbb{I}_l} \geq 0$. An algorithm overcoming this issue was proposed in~\cite{gill1986lssol} where the initial feasible point is determined by minimizing the sum of infeasibilities. However,  the algorithm fails when no feasible but only an optimal infeasible point  $A_{\mathbb{I}_l}x - b_{\mathbb{I}_l} < 0$ exists. The solvers proposed in~\cite{Kanoun2011,escande2014,dimitrov:2015} are based on the ASM and are able of handling all the above cases by virtue of relaxation with slack variables.
	Similarly, we introduce the notion of active constraints into~\eqref{eq:hlspLvll} (\eqref{eq:hlspLvll} is  equivalent to~\eqref{eq:hlsplexmin} only if all inequality constraints on all levels are feasible):
	\begin{align}
		\mini_{x,v_{\mathbb{E}_l},v_{\mathbb{I}_l}}& \qquad \frac{1}{2}\Vert v_{\mathbb{E}_l} \Vert^2 + \frac{1}{2}\Vert v_{\mathbb{I}_l} \Vert^2\qquad l=1,\dots,p \label{eq:hlsp}\tag{HLSP}\\
		\text{s.t.}
		& \qquad A_{\mathbb{E}_l}x - b_{\mathbb{E}_l} = v_{\mathbb{E}_l}\nonumber\\
		& \qquad A_{\mathbb{I}_l}x - b_{\mathbb{I}_l} \geq v_{\mathbb{I}_l}\nonumber\\
		& \qquad A_{\mA_{\cup l-1}}x - b_{\mA_{\cup l-1}} = v_{\mA_{\cup l-1}}^*\nonumber\\
		& \qquad A_{\mI_{\cup l-1}}x - b_{\mI_{\cup l-1}} \geq 0\nonumber
	\end{align}
	The active set $\mA_{\cup l-1}$ represents all equality constraints ${\bE_{\cup l-1}}$ and the inequality constraints of ${\bI_{\cup l-1}}$ that were infeasible and violated by $v_{\mA_{\cup l-1}}^*$ at the optimal points of the levels $1$ to $l-1$. Consequently, ${\mI_{\cup l-1}}$ are the remaining inequality constraints that were feasible and satisfied at the optimal points of the levels $1$ to $l-1$ and which are labeled as `inactive'. More details on setting these active sets are given in Sec.~\ref{sec:virtualPriorLvl}.
	
	The aim of the optimization problem is
	\begin{enumerate}
		\item to find the feasible or optimal infeasible point $v_{\mathbb{E}_l}^*$ and $v_{\mathbb{I}_l}^*$ of the equality and inequality constraints $\mathbb{E}_l$ and $\mathbb{I}_l$ of the current level $l$ \\
		$\rightarrow$ {IPM}
		\item while keeping the inactive inequality constraints $\mI_{\cup l-1}$ of the previous levels $1,\dots,l-1$ feasible with $A_{\mI_{\cup l-1}}x - b_{\mI_{\cup l-1}} \geq 0$ 	\\	$\rightarrow$ {IPM}
		\item and while keeping the active constraints $\mA_{\cup l-1}$ of the previous levels $1,\dots,l-1$ optimal at the current violation $v_{\mA_{\cup l-1}}^*$ \\
		$\rightarrow$ {Nullspace method}

	\end{enumerate}
	The first and second point we achieve by the IPM as outlined below.
	The third point is outlined in Sec.~\ref{sec:nsproj}.
	
	\subsection{The IPM for HLSP and the Newton's method}
	We introduce two positive slack variables, $w_{\mathbb{I}_l}$ for the inequality constraints on the current level $l$, and $w_{\mI_{\cup l-1}}$ for the inactive inequality constraints of the previous levels. They are then penalized by the log function to avoid values approaching zero (`log-barrier'):
	\begin{align}
		\min_x.& \qquad \frac{1}{2}\Vert v_{\mathbb{E}_l} \Vert^2 + \frac{1}{2}\Vert v_{\mathbb{I}_l} \Vert^2 - \sigma_{\mathbb{I}_l}\mu_{\mathbb{I}_l}\sum\log(w_{\mathbb{I}_l}) 
		- \sigma_{\mI_{\cup l-1}}\mu_{\mI_{\cup l-1}}\sum\log(w_{\mI_{\cup l-1}})\nonumber\\
		\text{s.t.}
		& \qquad A_{\mathbb{E}_l}x - b_{\mathbb{E}_l} = v_{\mathbb{E}_l}\nonumber\\
		& \qquad A_{\mathbb{I}_l}x - b_{\mathbb{I}_l} - v_{\mathbb{I}_l} = w_{\mathbb{I}_l}\nonumber\\
		&\qquad w_{\mathbb{I}_l} \geq 0 \nonumber\\
		& \qquad A_{\mA_{\cup l-1}}x - b_{\mA_{\cup l-1}} = v_{\mA_{\cup l-1}}^*\nonumber\\
		& \qquad A_{\mI_{\cup l-1}}x - b_{\mI_{\cup l-1}} = w_{\mI_{\cup l-1}}\nonumber\\
		&\qquad w_{\mI_{\cup l-1}} \geq 0\label{eq:ipmhqpopt} 
	\end{align}
	The Lagrangian of the optimization problem~\eqref{eq:ipmhqpopt} is
	\begin{align}
		\mathcal{L} &\coloneqq \frac{1}{2}\Vert v_{\mathbb{E}_l} \Vert^2 + \frac{1}{2}\Vert v_{\mathbb{I}_l} \Vert^2 - \sigma_{\mathbb{I}_l}\mu_{\mathbb{I}_l}\sum\log(w_{\mathbb{I}_l}) -\sigma_{\mI_{\cup l-1}}\mu_{\mI_{\cup l-1}}\sum\log(w_{\mI_{\cup l-1}})\nonumber\\
		& - \lambda_{\mathbb{E}_l}^T(A_{\mathbb{E}_l}x - b_{\mathbb{E}_l} - v_{\mathbb{E}_l}) - \lambda_{\mathbb{I}_l}^T(A_{\mathbb{I}_l}x - b_{\mathbb{I}_l} - v_{\mathbb{I}_l} - w_{\mathbb{I}_l})\nonumber\\
		&  - \lambda_{\mA_{\cup l-1}}^T (A_{\mA_{\cup l-1}}x - b_{\mA_{\cup l-1}} - v_{\mA_{\cup l-1}}^*)- \lambda_{\mI_{\cup l-1}}^T(A_{\mI_{\cup l-1}}x - b_{\mI_{\cup l-1}} - w_{\mI_{\cup l-1}})
	\end{align}
	$\lambda$ are the Lagrange multipliers associated with the corresponding constraints $\mathbb{E}_l$, $\mathbb{I}_l$, $\mA_{\cup l-1}$ and $\mI_{\cup l-1}$.
	
	The first order optimality or KKT conditions $K\coloneqq\nabla_q \mathcal{L} = 0$ with the variable vector $q$
	\begin{equation}
		q \coloneqq \BIN x^T & v_{\mathbb{E}_l}^T & v_{\mathbb{I}_l}^T & w_{\mathbb{I}_l}^T & \lambda_{\mA_{\cup l-1}}^T & \lambda_{\mI_{\cup l-1}} & w_{\mI_{\cup l-1}}^T \BOUT^T
	\end{equation}
	are (after applying the substitutions $v_{\mathbb{E}_l} = -\lambda_{\mathbb{E}_l}$ and $v_{\mathbb{I}_l} = -\lambda_{\mathbb{I}_l}$)
	\begin{align}
		K_l(q)&
		\hspace{-1pt}=	\hspace{-1pt}
		\BIN
		K_{x,l} \\
		K_{v_{\mathbb{E}_l},l} \\
		K_{v_{\mathbb{I}_l},l} \\
		K_{w_{\mathbb{I}_l},l} \\
		K_{\lambda_{\mA_{\cup l-1}},l} \\
		K_{\lambda_{\mI_{\cup l-1}},l} \\ 
		K_{w_{\mI_{\cup l-1}},l}
		\BOUT
		\hspace{-1pt}\coloneqq	\hspace{-1pt}
		\BIN
		A_{\mathbb{E}_l}^Tv_{\mathbb{E}_l} + A_{\mathbb{I}_l}^Tv_{\mathbb{I}_l} - A_{\mA_{\cup l-1}}^T\lambda_{\mA_{\cup l-1}} - A_{\mI_{\cup l-1}}^T\lambda_{\mI_{\cup l-1}}\\
		b_{\mathbb{E}_l} - A_{\mathbb{E}_l}x + v_{\mathbb{E}_l}  \\
		b_{\mathbb{I}_l} - A_{\mathbb{I}_l}x + v_{\mathbb{I}_l} + w_{\mathbb{I}_l}\\
		w_{\mathbb{I}_l}\odot v_{\mathbb{I}_l} + \sigma_{\mathbb{I}_l}\mu_{\mathbb{I}_l} e\\
		b_{\mA_{\cup l-1}} - A_{\mA_{\cup l-1}}x + v_{\mA_{\cup l-1}}^* \\
		b_{\mI_{\cup l-1}} - A_{\mI_{\cup l-1}}x + w_{\mI_{\cup l-1}}\\
		\lambda_{\mI_{\cup l-1}}\odot w_{\mI_{\cup l-1}} - \sigma_{\mI_{\cup l-1}}\mu_{\mI_{\cup l-1}} e
		\BOUT	\hspace{-1pt}=	\hspace{-1pt} 0
		\label{eq:kktHLSP}\tag{KKT}
	\end{align}
	The operator $\odot$ indicates the element-wise multiplication between two vectors.  $e$ is a vector of 1's of appropriate dimensions.
	The duality measures $\mu$ and the centering parameter $\sigma$ are given by~\citep{vanderbei2013}
	\begin{align}
		\mu_{\mathbb{I}_l} &\coloneqq \lambda_{\mathbb{I}_l}^Tw_i/(n+{m}_{l,i})\\
		\mu_{\mI_{\cup l-1}} &\coloneqq \lambda_{\mI_{\cup l-1}}^Tw_{\mI_{\cup l-1}} / (n + m_{\mI_{\cup l-1}})\\
		&\sigma_{\mathbb{I}_l}, \sigma_{m_{\mI_{\cup l-1}}}\in[0,1]
	\end{align}
	The values can also be determined by Mehtrotra's predictor-corrector algorithm~\citep{Mehrotra1992}, see Sec.~\ref{sec:predcoralgo}.
	
	The \ref{eq:kktHLSP} conditions are non-linear. We linearize them by the \text{Newton's method} by applying the Newton step
	\begin{equation}
		K_l(q + \Delta q) \approx K_l(q) + \nabla_q K_l(q)\Delta q
	\end{equation}
	The Lagrangian Hessian is given by
	\begin{align}
		&\nabla_q K_l(q) \coloneqq \BIN
		0 & A_{\mathbb{E}_l}^T & A_{\mathbb{I}_l}^T & 0 & -A_{\mA_{\cup l-1}}^T& -A_{\mI_{\cup l-1}}^T & 0\\
		-A_{\mathbb{E}_l} & I & 0 & 0 & 0 & 0 & 0\\ 
		-A_{\mathbb{I}_l} & 0 & I & I & 0 & 0 & 0\\
		0 & 0 & W_{\mathbb{I}_l} & V_{\mathbb{I}_l} & 0 & 0 & 0\\
		-A_{\mA_{\cup l-1}} & 0 & 0 & 0 & 0 & 0 & 0 \\
		-A_{\mI_{\cup l-1}} & 0 & 0 & 0 & 0 & 0 & I\\
		0 & 0 & 0 & 0 & 0 & W_{\mI_{\cup l-1}} & \Lambda_{\mI_{\cup l-1}}
		\BOUT
	\end{align}
	The capital variables $W = \diag(w) \in\mathbb{R}^{m,m}$, $V = \diag(V)\in\mathbb{R}^{m,m}$ and $\Lambda = \diag(\lambda)\in\mathbb{R}^{m,m}$ are the diagonal matrix equivalents of their vectors $w\in\mathbb{R}^{m}$, $v\in\mathbb{R}^{m}$ and $\lambda\in\mathbb{R}^{m}$. $I$ is an identity matrix of appropriate dimensions.
	
	We sequentially apply substitutions for $\Delta v_{\mathbb{E}_l}$, $\Delta v_{\mathbb{I}_l}$, $\Delta w_{\mI_{\cup l-1}}$, $\Delta w_{\mathbb{I}_l}$ and $\Delta \lambda_{\mI_{\cup l-1}}$
	which leads to the \textit{hierarchical augmented system}
	\begin{align}
		&\BIN C_l & -A_{\mA_{\cup l-1}}^T\\
		-A_{\mA_{\cup l-1}} & 0
		\BOUT
		\BIN
		\Delta x\\
		\Delta \lambda_{\mA_{\cup l-1}}
		\BOUT = \BIN r_{l,1} \\ r_{l,2} \BOUT 
		\label{eq:NeNmethodFull}
	\end{align}
	with
	\begin{align}
		C_l &\coloneqq A_{\mathbb{E}_l}^TA_{\mathbb{E}_l} + A_{\mathbb{I}_l}^T\left(I + ({V_{\mathbb{I}_l} - W_{\mathbb{I}_l}})^{-1}{W_{\mathbb{I}_l}}\right)A_{\mathbb{I}_l} + A_{\mI_{\cup l-1}}^Tw_{\mI_{\cup l-1}}^{-1}\lambda_{\mI_{\cup l-1}}A_{\mI_{\cup l-1}}
		\label{eq:cl}
	\end{align}
	and the right hand side
	\begin{align}
		r_{l,1} &\coloneqq A_{\mA_{\cup l-1}}^T\lambda_{\mA_{\cup l-1}} + A_{\mI_{\cup l-1}}^TF + A_{\mathbb{E}_l}^T(b_{\mathbb{E}_l} - A_{\mathbb{E}_l}x) + A_{\mathbb{I}_l}^TG 
		\\
		r_{l,2} &\coloneqq A_{\mA_{\cup l-1}}x - b_{\mA_{\cup l-1}} -w_{\mA_{\cup l-1}}^*\label{eq:r2}
	\end{align}
	$G$ and $F$ are given by
	\begin{align}
		F_l &\coloneqq \lambda_{\mI_{\cup l-1}}
		+ W_{\mI_{\cup l-1}}^{-1}(\lambda_{\mI_{\cup l-1}}\odot(b_{\mI_{\cup l-1}}- A_{\mI_{\cup l-1}}x) \label{eq:F}+ \sigma_{\mI_{\cup l-1}}\mu_{\mI_{\cup l-1}} e)\\
		G_l &\coloneqq b_{\mathbb{I}_l} - A_{\mathbb{I}_l}x + w_{\mathbb{I}_l} \label{eq:G}- ({V_{\mathbb{I}_l} - W_{\mathbb{I}_l}})^{-1}({\sigma_{\mathbb{I}_l}\mu_{\mathbb{I}_l} e + w_{\mathbb{I}_l}\odot(A_{\mathbb{I}_l}x - b_{\mathbb{I}_l} - w_{\mathbb{I}_l})})
	\end{align}
	respectively.
	
	\subsection{The iterative nature of the IPM}
	The linear hierarchical augmented system~\eqref{eq:NeNmethodFull} of level $l$ is now repeatedly solved for the primal and dual steps $\Delta q$.
	The IPM's for finding the initial feasible or optimal infeasible point require following \textit{dual feasibility} conditions to hold
	\begin{equation}
		v_{\mathbb{I}_l} \leq 0 \qquad \text{and} \qquad
		w_{\mathbb{I}_l} \geq 0\label{eq:linesearch_li}
	\end{equation}
	For the IPM of the inactive constraints we have
	\begin{equation}
		\lambda_{\mI_{\cup l-1}} \geq 0 \qquad \text{and} \qquad
		w_{\mI_{\cup l-1}} \geq 0\label{eq:linesearch_l-1i}
	\end{equation}
	These conditions are maintained by line search. 
	The computed primal and dual steps  $\Delta q$ are then scaled by the line search factor $\alpha$ and added to the current estimates of the primal and dual $q\leftarrow q + \alpha \Delta q$. The Newton's method terminates once the norm of the non-linear \ref{eq:kktHLSP} conditions $\left\Vert K_l(q)\right\Vert_2 < \epsilon$~\eqref{eq:kktHLSP} is below a certain threshold $\epsilon = 10^{-12}$.
	
	\section{Computing the Newton step for IPM for HLSP}
	\label{sec:compHLSP}
	
	Solving the augmented system~\eqref{eq:NeNmethodFull} directly in a dense manner is inefficient as the zero block in the lower right corner of the left hand side matrix would be ignored (a sparse solver has been proposed for example in~\cite{qpSwift2019}).  One way to circumvent this is to form the Schur complement~\citep{Bartlett2006,zeilinger2012}
	\begin{align}
		\Delta x &= C_l^{-1}(r_{l,1} + A_{\mA_{\cup l-1}}^T \Delta \lambda_{\mA_{\cup l-1}})
		\label{eq:schurdx}
	\end{align}
	The resulting `Schur' normal equations are given below (IPM-snf: Schur normal form or Schur normal equations):
	\begin{equation}
		A_{\mA_{\cup l-1}}
		C_l^{-1}A_{\mA_{\cup l-1}}^T 
		\Delta \lambda_{\mA_{\cup l-1}}
		=
		-r_{l,2} - A_{\mA_{\cup l-1}}C_l^{-1}r_{l,1}
		\tag{IPM-snf}
		\label{eq:qpclassicalNeq}
	\end{equation}
	In order to obtain the primal step $\Delta x$ two decompositions per Newton iteration are required, one for $C_l^{-1}$ and one for $A_{\mA_{\cup l-1}}C_l^{-1}A_{\mA_{\cup l-1}}^T$.  Efficient strategies especially in the context of Model-Predictive-Control (MPC) leveraging banded matrix structures have been proposed~\citep{zeilinger2012}. However, this still can prove inefficient, especially if the dual $\lambda_{\mA_{\cup l-1}}$ is not required for later use (they might be required for example for Hessian computations~\citep{pfeiffer2023}). We therefore describe in the following the algorithm $\mathcal{N}$IPM-HLSP based on the nullspace method which only requires a single decomposition per Newton iteration. 
	
	\subsection{Algorithm for $\mathcal{N}$IPM-HLSP}
	\label{sec:algo}

	An overview of our algorithm $\mathcal{N}$IPM-HLSP to resolve~\ref{eq:hlsp}'s is given below:
	\begin{enumerate}
		\item Go through the hierarchy~\ref{eq:hlsp} from level $1$ to $p$.
		\item For each level $l$, repeatedly compute the Newton step by solving~\ref{eq:HLSPNeNmethod} or~\ref{eq:LQNmethod} formulated in Sec.~\ref{sec:nsproj} and Sec.~\ref{sec:lqipm}, respectively. In case of the Mehrotra's predictor corrector algorithm (Sec.~\ref{sec:predcoralgo}), first compute the affine step $\Delta \_^{\text{aff}}$ and then the centered one $\Delta \_$. The symbol $\_$ is a placeholder for the different variables $x, v_{\mathbb{I}_l}, w_{\mathbb{I}_l},w_{\mI_{\cup l-1}}$ and $\lambda_{\mI_{\cup l-1}}$.
		\item Line search for dual feasibility~\eqref{eq:linesearch_li} and~\eqref{eq:linesearch_l-1i}. Add step scaled by the line search factor to current primal and dual estimate.
		\item Upon convergence of the Newton's method $\Vert \tilde{K} \Vert \leq \epsilon$, gather all inactive inequality constraints from higher priority levels ${\mI_{\cup l-1}}$ that are saturated and add them to the active set $\mathcal{A}_{l^*}$ of the virtual priority level $l^*$. The active set assembly including the concept of virtual priority levels is explained in Sec.~\ref{sec:virtualPriorLvl}. $\tilde{K}$ are the projected \ref{eq:kktHLSP} conditions described in Sec.~\ref{sec:ipmnodual}. 
		\item Compute its nullspace $Z_{\mA_{l^*}},r\leftarrow\mathcal{N}(A_{\mA_{l^*}})$ (the operator $\mathcal{N}(A)$ returns the rank and a basis of the nullspace $Z$ of the input matrix $A$) and project lower priority levels and the remaining inactive inequality constraints into it. Augment $N_{\mA_{\cup l^*}} \leftarrow N_{\mA_{\cup l-1}}Z_{\mathcal{A}_{l^*}}$.
		\item Add all equalities and the violated inequality constraints from level $l$ to the active set $\mathcal{A}_l$ of level $l$. This is described in Sec.~\ref{sec:virtualPriorLvl}.
		\item Compute its nullspace $Z_{\mathcal{A}_{l}},r\leftarrow\mathcal{N}(A_{\mathcal{A}_{l}})$ and project lower priority levels and the remaining inactive inequality constraints into it. Augment $N_{\mA_{\cup l}} \leftarrow N_{\mA_{\cup l^*}}Z_{\mathcal{A}_{l}}$.
	\end{enumerate}
	A pseudo-implementation of the algorithm is given in App.~\ref{app:alg}.

	\subsection{$\mathcal{N}$IPM-HSLP: The nullspace method based IPM for HLSP}
	\label{sec:nsproj}

	We first assume that $r_{l,2} = 0$~\eqref{eq:r2}.
	This condition holds as $x$ is updated after every Newton iteration during the resolution of the higher priority levels $l=1,\dots,l-1$ and is therefore feasible with respect to the active constraints $\mA_{\cup l-1}$. 
	We then apply the nullspace method~\citep{Nocedal2006} by first introducing the variable change
	\begin{equation}
		\Delta x = N_{\mA_{\cup l-1}}\Delta z
		\label{eq:HLSPchangeOfVar}
	\end{equation}
	The nullspace basis $N_{\mA_{\cup l-1}}$ of $A_{\mA_{\cup l-1}}$ fulfills the condition $A_{\mA_{\cup l-1}}N_{\mA_{\cup l-1}} = 0$. This means that $A_{\mA_{\cup l-1}}\Delta x = 0$ such that the condition $r_{l,2} = 0$ continues to be fulfilled.
	An additional multiplication of the first row of~\eqref{eq:NeNmethodFull} by $N_{\mA_{\cup l-1}}^T$ from the left results in the `projected' normal equations
	\begin{align}
		N_{\mA_{\cup l-1}}^TC_lN_{\mA_{\cup l-1}} 
		\Delta z
		=
		N_{\mA_{\cup l-1}}^Tr_{l,1}
		\tag{$\mathcal{N}$IPM-nf}
		\label{eq:HLSPNeNmethod}
	\end{align}
	As can be seen, the primal step $\Delta z$ is directly obtained without the intermediate step of computing the Lagrange multipliers $\lambda_{\mA_{\cup l-1}}$ corresponding to the active constraints $\mA_{\cup l-1}$. They can be obtained by solving
	\begin{equation}
		A_{\mA_{\cup l-1}}^T\Delta \lambda_{\mA_{\cup l-1}}= 
		C_l\Delta x - r_{l,1} 
		\label{eq:lambdaact}
	\end{equation}
	An efficient recursive method of calculation is detailed in Appendix~\ref{app:recLambdaact}. As we detail in Sec.~\ref{sec:ipmnodual}, we do not necessarily need to compute these Lagrange multipliers after all.
	
	We use the nullspace basis described in~\cite{Nocedal2006}
	\begin{align}
		Z_{\mathcal{A}_{l}} = 
		P
		\BIN
		-R^{-1}T\\
		I
		\BOUT
		\quad 
		\text{with}\quad
		A_{\mathcal{A}_{l}} = 
		Q
		\BIN
		R & T\\
		0 & 0
		\BOUT
		P^T	
		\label{eq:nsb}
	\end{align}
	Resulting from a rank-revealing QR decomposition (RRQR) of $A_{\mathcal{A}_{l}}$, $Q$ is an orthogonal matrix, $R$ is upper triangular, $T$ is a rectangular matrix and $P$ is a permutation matrix. The bottom zero row is due to possible linear dependencies in $A_{\mathcal{A}_l}$. 
	This basis has variable reducing qualities~\citep{bjoerck1996} (i.e. projected matrices are reduced in variables), projections are computed cheaply due to the lower identity matrix and it can be reused for the calculation of the Lagrange multipliers (App.~\ref{app:recLambdaact}). While this basis is partly sparse, its projections are dense if there is any density in the factors $R$ or $T$. 	The accumulated nullspace basis $N_{\mathcal{A}_{\cup l}}$ is computed by $N_{\mathcal{A}_{\cup l}} = Z_{\mathcal{A}_{1}} \dots Z_{\mathcal{A}_{l}}$.

	\subsection{Active set assembly and virtual priority levels}
	\label{sec:virtualPriorLvl}
	
	Once the Newton's method of a priority level $l$ has converged, inactive constraints $\mI_{\cup l-1}$ from previous levels $1$ to $l-1$ are either saturated with $A_{\mI_{\cup l-1}}x - b_{\mI_{\cup l-1}} = 0$ or satisfied with $A_{\mI_{\cup l-1}}x - b_{\mI_{\cup l-1}} > 0$. Additionally, inequality constraints $\mathbb{I}_l$ from the current level $l$ are either satisfied with $A_{\mathbb{I}_l}x - b_{\mathbb{I}_l} \geq 0$ or violated with $A_{\mathbb{I}_l}x - b_{\mathbb{I}_l} < 0$. In order to be able to proceed with the cascade-like resolution of the next priority level we need to determine two distinct active sets $\mathcal{A}_{l^*}$ and $\mathcal{A}_{l}$ (such that $\mA_{\cup l} \coloneqq  \mathcal{A}_{1^*} \cup \mathcal{A}_1 \cup \cdots \cup \mathcal{A}_{l^*} \cup \mathcal{A}_l $).
	The index $l^*$ represents a `virtual' priority level whose active set $\mathcal{A}_{l^*}$ contains saturated constraints from $\mI_{\cup l-1}$. This is necessary as the projection of $\mI_{\cup l-1}$ into the nullspace basis of the active constraints $\mathcal{A}_l$ of level $l$ would contradict the priority order. Therefore, a virtual priority level $l^*$ has lower priority than level $l-1$ but higher priority than level $l$ whose active set $\mA_l$ comprises of the equality constraints $\bE_l$ and the activated constraints from $\bI_l$. 
	
	We activate a constraint $c \in\mI_{\cup l-1}$ if the corresponding pair of slack and Lagrange multiplier fulfills the conditions
	\begin{align}
		w_{\mI_{\cup l-1}}(c) < \xi \qquad \text{ and } \qquad \lambda_{\mI_{\cup l-1}}(c) > \xi 
	\end{align}
	While $w_{\mI_{\cup l-1}}(c) < \xi$ ensures that only saturated constraints are activated, the condition $\lambda_{\mI_{\cup l-1}}(c) > \xi $ only selects constraints that are in significant conflict with constraints from level $l$ and are not just accidentally saturated, for example by a randomly chosen initial point. Otherwise the ability of resolving the lower priority level $l$ is artificially restricted by an unnecessarily inflated active set $\mA_{\cup l-1}$.
	
	The threshold $\xi = 10^{-8}$ is necessary as the Newton's method is usually not run to complete convergence $\Vert K_l\Vert_2 = 0 $ but stopped earlier with some threshold $\Vert K_l\Vert_2 < \epsilon$. The small choice for $\xi$ requires the Newton's method to converge to a similarly accurate degree $\epsilon = 10^{-12}$ in order to have an accurate value of the dual $w_{\mI_{\cup l-1}}$ and $\lambda_{\mI_{\cup l-1}}$. 
	
	The following steps are conducted for the virtual priority level $l^*$:
	\begin{enumerate}
		\item Add the newly activated constraints from ${\mI_{\cup l-1}}$ to  ${\mathcal{A}_{l^*}}$ as a `virtual' priority level $l^*$.
		\item Calculate the nullspace basis $Z_{\mA_{l^*}}$ of $A_{\mathcal{A}_{l^*}}N_{\mA_{\cup l-1}}$ and project lower priority levels and the remaining inactive inequality constraints into it.
		\item Augment $N_{\mA_{\cup l-1}}$ to $N_{\mA_{\cup l^*}} = N_{\mA_{\cup l-1}}Z_{\mA_{l^*}}$
	\end{enumerate}
	Now we handle the violated constraints from level $l$ for the active set $\mathcal{A}_{l}$. A constraint $c\in\mathbb{I}_l$ is activated if the following conditions are fulfilled:
	\begin{equation}
		w_{\mathbb{I}_l}(c) < \xi \qquad \text{ and } \qquad v_{\mathbb{I}_l}(c) < -\xi 
	\end{equation}
	Following steps complete the active set assembly:
	\begin{enumerate}
		\item Assemble the active set $\mathcal{A}_l$ of level $l$ with all equality constraints and the active inequality constraints.
		\item Inactive constraints are added to ${\mI_{\cup l}}$. Note that in the case of bound constraints we check whether the variable is already constrained by the same constraint. If this is the case, no new constraint is added and the tighter bound of the two is used as the new right hand side $b_{\mI_{\cup l}}$. This prevents unnecessarily increasing the size of ${\mI_{\cup l}}$. Note that this is already implied by the set union symbol $\cup$.
		\item Calculate the nullspace basis $Z_{\mA{l}}$ of $A_{\mathcal{A}_{l}}N_{\mA_{\cup l^*}}$ and project lower priority levels and the remaining inactive inequality constraints into it.
		\item Augment $N_{\cup\mA_{l^*}}$ to $N_{\mA_{\cup l}} = N_{\mA_{\cup  l^*}}Z_{\mA_l}$.
	\end{enumerate}

	\subsection{Avoiding the calculation of $\Delta \lambda_{\mA_{\cup l-1}}$}
	\label{sec:ipmnodual}
	
	The cost of calculating the dual step $\Delta \lambda_{\mA_{\cup l-1}}$ associated with the equality constraints ${\mA_{\cup l-1}}$ is of magnitude $O(r_{\mA_{\cup l-1}}^2)$. The rank of $A_{\mA_{\cup l-1}}$ is given as $r_{\mA_{\cup l-1}}$. However, we do \textit{not necessarily} need to update the Lagrange multipliers $\lambda_{\mA_{\cup l-1}}$ since none of the other primal or dual variables depend on them. Additionally, we can explicitly calculate $\lambda_{\mA_{\cup l-1}}$ by solving $K_{x,l} = 0$~\eqref{eq:kktHLSP}. This is in contrast to the augmented system~\eqref{eq:NeNmethodFull} or the Schur normal equations~\eqref{eq:qpclassicalNeq} which require the calculation of the dual step in order to obtain the primal step. 
	
	The dual $\lambda_{\mA_{\cup l-1}}$ is only necessary for the evaluation of the norm of $K_l$~\eqref{eq:kktHLSP} to determine whether the IPM has converged with $\Vert K_{l} \Vert_2 < \epsilon$. However, we can instead use the projected \ref{eq:kktHLSP} conditions $\tilde{K}$ (we only show the relevant components)
	\begin{align}
		\tilde{K}_{x,l} &\coloneqq N_{\mA_{\cup l-1}}^T(A_{\mathbb{E}_l}^Tv_{\mathbb{E}_l} + A_{\mathbb{I}_l}^Tv_{\mathbb{I}_l} - A_{\mI_{\cup l-1}}^T\lambda_{\mI_{\cup l-1}})\\
		\tilde{K}_{\lambda_{\mA_{\cup l-1}},l}& \coloneqq	N_{\mA_{\cup l-1}}^T(b_{\mA_{\cup l-1}} - A_{\mA_{\cup l-1}}x + v_{\mA_{\cup l-1}}^*) = 0
	\end{align}
	The nullity of the second equation holds since we already have ensured primal feasibility $v_{\mA_{\cup l-1}}^*$ of the active constraints $\mA_{\cup l-1}$ when resolving the previous levels $1,\dots,l-1$. Furthermore, any nullspace step $\Delta x = N_{\mA_{\cup l-1}}\Delta z$ does not influence ${K}_{\lambda_{\mA_{\cup l-1}},l}$ since $A_{\mA_{\cup l-1}}(x + N_{\mA_{\cup l-1}}\Delta z) = A_{\mA_{\cup l-1}}x$.

	\subsection{The least-squares form of the $\mathcal{N}$IPM-HLSP}
	\label{sec:lqipm}
	
	\ref{eq:HLSPNeNmethod} can be rewritten to least squares form
	\begin{align}
		\mini_{\Delta z}\left\Vert\BIN
		\sqrt{w_{\mI_{\cup l-1}}^{-1}\lambda_{\mI_{\cup l-1}}}\utA_{\mI_{\cup l-1}}\\
		\sqrt{I + ({V_{\mathbb{I}_l} - W_{\mathbb{I}_l}})^{-1}{W_{\mathbb{I}_l}}}\tA_{\mathbb{I}_l}\\
		\tA_{\mathbb{E}_l}
		\BOUT
		\hspace{-2pt}
		\Delta z
		\hspace{-2pt}-\hspace{-2pt}
		\BIN
		\sqrt{w_{\mI_{\cup l-1}}\lambda_{\mI_{\cup l-1}}^{-1}}F\\
		\sqrt{I + ({V_{\mathbb{I}_l} - W_{\mathbb{I}_l}})^{-1}{W_{\mathbb{I}_l}}}^{-1}G\\
		b_{\mathbb{E}_l}-A_{\mathbb{E}_l}x
		\BOUT
		\right\Vert^2
		\tag{$\mathcal{N}$IPM-ls}
		\label{eq:LQNmethod}
	\end{align}
	$F$ and $G$ are defined in~\eqref{eq:F} and~\eqref{eq:G}, respectively. We use the notation $\tilde{M} \coloneqq MN_{\mA_{\cup l-1}}$. The above is well defined as we show in the following.
	
	\begin{theorem}
		The expressions under the square root of~\ref{eq:LQNmethod} are non-negative.
	\end{theorem}
	\begin{proof}
		$W_{\mI_{\cup l-1}}^{-1}\Lambda_{\mI_{\cup l-1}}\geq0$ and $W_{\mI_{\cup l-1}}\Lambda_{\mI_{\cup l-1}}^{-1}\geq0$ follow directly from~\eqref{eq:linesearch_l-1i}.
		Furthermore, for $\sqrt{I + ({V_{\mathbb{I}_l} - W_{\mathbb{I}_l}})^{-1}{W_{\mathbb{I}_l}}}$ we have ${V_{\mathbb{I}_l} - W_{\mathbb{I}_l}}\leq 0$ because of the dual feasibility conditions~\eqref{eq:linesearch_li}. This leads to 
		\begin{align}
			I + ({V_{\mathbb{I}_l} - W_{\mathbb{I}_l}})^{-1}{W_{\mathbb{I}_l}} \geq 0 \quad\leftrightarrow\quad
			V_{\mathbb{I}_l} \leq 0
		\end{align}
		which is true and ensured by~\eqref{eq:linesearch_li}.
	\end{proof}
	
	The least-squares form has the advantage that the (substituted) Lagrangian Hessian~\eqref{eq:cl} does not need to be computed. On the other hand, solving the system~\ref{eq:LQNmethod} requires a more expensive rank revealing decomposition (for example QR decomposition). 
	Following our considerations from Sec.~\ref{sec:qpopcomp}, we use the condition 
	\begin{equation}
		n_r > \frac{3}{5} (m_{\mI_{\cup l-1}} + m_{\mathbb{E}_l} + m_{\mathbb{I}_l})
		\label{eq:LScond}
	\end{equation}
	to decide whether the least-squares form~\eqref{eq:LQNmethod} is more appropriate than the projected normal equations~\eqref{eq:HLSPNeNmethod}. 
	
	The above implies that the least-squares form is more efficient for a high ratio of number of variables to number of constraints ($n_r \gg m$). Yet, efficient parallel decomposition methods for tall matrices ($m \gg n_r$) have been evaluated for example in~\cite{Sahinidis2020}.
	Furthermore, the above condition neglects the circumstance that the least-squares form can potentially maintain a high degree of sparsity of the constraint matrices $A$. In future work we therefore aim to implement a sparse version of the solver, for example based on sparsity maintaining nullspace bases as proposed in~\cite{gill1987}.

	\subsection{Mehrotra's predictor corrector algorithm}
	\label{sec:predcoralgo}
	
	In the following we adapt Mehrotra's predictor-corrector algorithm~\citep{Mehrotra1992} to~\ref{eq:hlsp}. It enables fast convergence of the Newton's method by choosing $\sigma_{\mI_{\cup l-1}}$ and $\sigma_{\mathbb{I}_l}$ appropriately. The algorithm for LSP is described in~\cite{Nocedal2006} (Alg. 16.4) and only requires slight adaptations for HLSP which are detailed in Appendix~\ref{app:predcor}.
	
	Even though the primal ($O(n_r^2)$) needs to be computed twice per Newton iteration this effort is negligible in comparison to the cost of computing a decomposition ($O(n_r^3)$) of the linear systems~\ref{eq:HLSPNeNmethod} or~\ref{eq:LQNmethod}.

	\section{Computational complexity of $\mathcal{N}$IPM-HLSP}
	\label{sec:qpopcomp}
	
	\begin{table}[htp!]
		\centering
		\caption{Necessary steps (\#) per Newton iteration and their number of operations  for the different IPM formulations for HLSP.}
		\begin{tabular}{@{} m{0.015\columnwidth}m{0.29\columnwidth}m{0.29\columnwidth}m{0.29\columnwidth} @{}}  
			\toprule
			\# & \ref{eq:qpclassicalNeq} & \ref{eq:HLSPNeNmethod} & \ref{eq:LQNmethod}\\
			\toprule
			1. & Form  $C_l$ \eqref{eq:cl}:
			
			$O((m_{\mathbb{E}_l} + m_{\mathbb{I}_l} + m_{\mI_{\cup l-1}})n^2)$ & \cellcolor{gray!10} Form $N_{\mA_{\cup l-1}}^T C_l N_{\mA_{\cup l-1}}$:
			
			$O((m_{\mathbb{E}_l} + m_{\mathbb{I}_l} + m_{\mI_{\cup l-1}})n_r^2)$ & QR decomposition of~\ref{eq:LQNmethod}: 
			
			$O(2n_r^2(m_{\mI_{\cup l-1}} + m_{\mathbb{E}_l} + m_{\mathbb{I}_l}) - 2n_r^3/3)$\\
			2. & \cellcolor{gray!10}LDLT decomposition of $C_l$:
			$O(n^3/3)$ & LDLT decomposition of $N_{\mA_{\cup l-1}}^T C_l N_{\mA_{\cup l-1}}$: $O(n_r^3/3)$ & \cellcolor{gray!10}Solve for $\Delta z$:  $O(n_r^2)$\\
			\ 3. & Form matrix $A_{\mA_{\cup l-1}}
			C_l^{-1}A_{\mA_{\cup l-1}}^T$:
			
			$O(m_{\mA_{\cup l-1}}n^2 + m_{\mA_{\cup l-1}}^2n)$ \citep{zeilinger2012} & \cellcolor{gray!10}Solve for $\Delta z$: $O(n_r^2)$ &  $\Delta x = N\Delta z$: $O(nn_r)$\\
			4. & \cellcolor{gray!10}LDLT decomposition of $A_{\mA_{\cup l-1}}
			C_l^{-1}A_{\mA_{\cup l-1}}^T$:
			
			$O(m_{\mA_{\cup l-1}}^3/3)$ &  $\Delta x = N_{\mA_{\cup l-1}}\Delta z$: $O(nn_r)$& \cellcolor{gray!10}\\
			5. & Calculate dual step $\Delta \lambda_{\mA_{\cup l-1}}$ \eqref{eq:qpclassicalNeq}:
			
			$O(m_{\mA_{\cup l-1}}^2)$ &\cellcolor{gray!10} & \\
			6. & \cellcolor{gray!10}Calculate primal step $\Delta x$ \eqref{eq:schurdx}: $O(n^2)$ & & \cellcolor{gray!10}\\
			\bottomrule
			$\sum$ & $O((m_{\mathbb{E}_l} + m_{\mathbb{I}_l} + m_{\mI_{\cup l-1}} + n/3 + m_{\mA_{\cup l-1}} + 1)n^2 + (n + m_{\mA_{\cup l-1}}/3 + 1)m_{\mA_{\cup l-1}}^2)$ &\cellcolor{gray!10} $O((m_{\mathbb{E}_l} + m_{\mathbb{I}_l} + m_{\mI_{\cup l-1}} + 1 + n_r)n_r^2 + nn_r)$ & $O(2n_r^2(m_{\mI_{\cup l-1}} + m_{\mathbb{E}_l} + m_{\mathbb{I}_l}) - 2n_r^3/3 + n_r^2 + nn_r)$\\
			\bottomrule
		\end{tabular}
		\label{tab:NewtonIterOp}
	\end{table}

	\subsection{Computational complexity of single Newton iteration}
	
	Table~\ref{tab:NewtonIterOp} details the number of operations per Newton iteration for the Schur~\eqref{eq:qpclassicalNeq} and projected normal equations~\eqref{eq:HLSPNeNmethod} and the least squares form~\eqref{eq:LQNmethod}.
	The projected normal equations (4 steps) and the least squares form (3 steps) require less steps per Newton iteration than the Schur normal equations (6 steps).
	However, both the {projected normal equations}~\eqref{eq:HLSPNeNmethod} and the least squares form~\eqref{eq:LQNmethod} require the calculation of a nullspace basis $Z_{\mA_{l-1}}$ of $A_{\mA_{\cup l-1}}$ in $O(2m_{\mA_{\cup l-1}}^2n - 2/3m_{\mA_{\cup l-1}}^3)$  at the beginning of the Newton's method of each level. Furthermore, the projections $A_{\mathbb{E}_l}Z_{\mA_{ l-1}}$ ($O(2m_{\mathbb{E}_l}nn_r)$), $A_{\mathbb{I}_l}Z_{\mA_{ l-1}}$ ($O(2m_{\mathbb{I}_l}nn_r)$) and $A_{\mI_{\cup l-1}}Z_{\mA_{ l-1}}$ ($O(2m_{\mI_{\cup l-1}}nn_r)$) are required.
	Nonetheless, with the additional effect of variable reduction due to the nullspace projections~\citep{bjoerck1996} we demonstrate in our evaluation (Sec.~\ref{eval:eqonly}) that $\mathcal{N}$IPM-HLSP is equivalently fast or faster than the Schur normal form. This holds even for equality only problems where convergence is achieved within one Newton iteration without offsetting the overhead of the nullspace method over several Newton iterations.
	
	\subsection{$\mathcal{N}$IPM-HLSP's overall performance}	
	
	A measure for the operational cost of $\mathcal{N}$IPM-HLSP dependent on the number of Newton iterations is given by
	\begin{align}
		&c_{\mathcal{N}\text{IPM-HLSP}} \coloneqq\\ 
		&\sum_{l=1}^{p}
		\iota_l \cdot c_{\text{IPM}}(n_r, m_{\mI_{\cup l-1}} + m_{\mathbb{I}_l}, m_{\mathbb{E}_l})
		+ c_{\mathcal{N},\text{IPM}}(m_{\mathcal{A}_{l^*}},m_{\mathcal{A}_{l}},m_{\mI_{\cup l}},m_{l+1:p}, n_r)\nonumber
	\end{align} 
	The first component reflects the cost $c_{\text{IPM}}$ of resolving the Newton's method of each level with $\iota_l$ Newton iterations  as is listed in Tab.~\ref{tab:NewtonIterOp}. The main computational load originates from solving the linear systems~\ref{eq:HLSPNeNmethod} or~\ref{eq:LQNmethod}. The cost is a function of the number of remaining variables $n_r$, the number of active constraints $m_{\mA_{l^*}}$ and $m_{\mA_{l}}$ of each respective level, the number of inactive constraints from the previous levels $m_{\mI_{\cup l-1}}$, the number of equalities and inequalities of the current level $m_{\mathbb{I}_l}$ and $m_{\mathbb{E}_l}$ and finally the number of all equality and inequality constraints from the remaining levels $m_{l+1:p}$. Secondly, once the Newton's method of level $l$ converges, the nullspace bases of the active-sets $\mathcal{A}_{l^*}$ and $\mathcal{A}_{l}$  have to be computed. Furthermore, the remaining levels $l+1,\dots,p$ and inactive constraints $\mI_{\cup l}$ need to be projected into it. This is reflected in the nullspace projection cost $c_{\mathcal{N},\text{IPM}}$.

	In comparison, the pseudo-cost of the ASM in resolving the HLSP is composed as follows:
	\begin{align}
		c_{\text{ASM}} & \coloneqq \sum_{i=1}^{\iota}\sum_{l=1}^{p} c_{\mathcal{N},\text{ASM}}(m_{\mathcal{A}_{l}},m_{l+1:p}, n_r)
	\end{align}
	$\iota$ are the number of active set iterations. Note that the cost of the ASM $c_{\text{ASM}}$ is included in the projection cost $c_{\mathcal{N},\text{ASM}}$ (computation of a rank revealing QR decomposition of the current level with ${m}_{\mathcal{A}_{l}}$ active constraints; ${m}_{\mathcal{A}_{l}}$ also counts the activated constraints from higher priority levels 1 to $l-1$). As can be seen, the nullspace projection of the whole active set needs to be repeated for each active set iteration whereas $\mathcal{N}$IPM-HLSP needs to project the whole HLSP (including inactive constraints) only once per HLSP resolution. This leads to less operations for $\mathcal{N}$IPM-HLSP in case of large number of active set iterations for the ASM. Note that unlike the IPM, the ASM does not need to consider the $m_{\bI_{\cup l}}$ inactive inequality constraint. This can be computationally advantageous depending on the problem's dimensions (see condition~\eqref{eq:LScond} as a reference).
	
	\section{Evaluation}
	\label{sec:eval}
	
	We employ our solver $\mathcal{N}$IPM-HLSP (\ref{eq:HLSPNeNmethod}: projected normal form,~\ref{eq:LQNmethod}: least-squares form) first on equality only constrained HLSP's of two levels $p=2$ to gain some insights into its computational efficiency (Sec.~\ref{eval:eqonly}). Especially, we are interested in how it fares in comparison to the Schur normal equations (\ref{eq:qpclassicalNeq}: Schur normal form). 
	
	Secondly, we use the solver within the sequential hierarchical least-squares programming solver with trust region (S-HLSP) presented in~\cite{pfeiffer2023} to solve NL-HLSP's of the form (as a sub-form of~\ref{eq:HNLP})
	\begin{align}
		\min_{u,{v}_{\bE_l}{v}_{\bI_l}} \quad & \frac{1}{2} \left\|{v}_{\bE_l}\right\|^2 + \frac{1}{2} \left\|{v}_{\bI_l}\right\|^2 \qquad\quad\hspace{8pt} l=1,\dots ,p\label{eq:NL-HLSP}\tag{NL-HLSP}
		\\
		\mbox{s.t.} \quad 
		& {f}_{\bE_l}({u}) \hspace{3pt} = \hspace{3pt} {v}_{\bE_l} \nonumber\\
		& {f}_{\bI_l}({u}) \hspace{3pt} \geq \hspace{3pt} {v}_{\bI_l} \nonumber\\
		& {f}_{\mA_{\cup l-1}}  =  {v}_{\mA_{\cup l-1}}^* \nonumber\\
		& {f}_{\mI_{\cup l-1}}  \geq  0 \nonumber
	\end{align}
	$f\in\mathbb{R}^m$ are non-linear task functions dependent on the variable vector $u\in\mathbb{R}^n$. $v\in\mathbb{R}^m$ are slack variables and $v^*\in\mathbb{R}^m$ are the optimal ones.
	S-HLSP repeatedly linearizes the \ref{eq:NL-HLSP} to \ref{eq:hlsp} `sub-problems' around the current working point $u$ (with the constraint matrices $A$ and vectors $b$ representing this linearization, for example as Jacobians and Hessians of $f$). The resulting step $x$ from a~\ref{eq:hlsp} sub-problem solver (for example $\mathcal{N}$IPM-HLSP) is then used to make an approximate step $u \leftarrow u + x$ within the non-linear model. The trust region radius limits the step $\Vert x\Vert_{\infty} < \rho$ in order to maintain validity of the linear approximation~\eqref{eq:hlsp} with respect to the original problem~\eqref{eq:NL-HLSP}. In all following examples, the trust region radius constraint is formulated on level 0 of the HLSP's. S-HLSP employs a real-time capable trust-region radius adaptation method with local convergence properties. Other non-linear programming methods like filter methods with global convergence properties~\citep{fletcher2002b} require the recomputation of sub-problem steps $x$ with different trust-region radii. Since this is not possible in real-time control due to computational limitations, every step is accepted. Instead, the trust region radius is only adapted (increased for valid step, decreased for invalid step) in the following control iteration (for more details see~\cite{pfeiffer2018}). 
	
	Given the conceptual context of S-HLSP as in~\cite{pfeiffer2023}, we accordingly sample our simulations from different robot control scenarios. All constraints $f(u)$ therefore consist of either non-linear trigonometric functions or linear bound constraints on the robot variables $u$ (joint velocities, joint torques, contact forces). This lets us design ill-posed problems commonly found in instantaneous robot control resulting from singular kinematic configurations of the robot. In these situations, first-order linearizations (Jacobians of the constraints $f(u)$) of the~\ref{eq:NL-HLSP} become ill-posed or singular. Effective solutions to counter these situtations have been proposed, for example in form of regularization~\citep{Chiaverini1997} or higher-order approximations~\citep{pfeiffer2023}.  For all robot control scenarios, one~\ref{eq:hlsp} sub-problem is solved per control iteration.
	
	We first evaluate a kinematic robot control example ($p=5$, see Sec.~\ref{sec:eval:clasproj}). 
	The evaluation is continued with a dynamic robot simulation depicting the effect of a larger number of variables and inequality constraints and ill-posed constraints due to singularities on the solver timings ($p=6$, Sec.~\ref{sec:eval:nonsing}).
	The evaluation is concluded with a push simulation in order to observe the distinctive difference in behavior between the ASM and IPM in case of a singular instance of a large change of the active set ($p=5$, Sec.~\ref{sec:eval:dynPush}).
	For the simulations we use the HRP-2Kai humanoid robot with 38 degrees of freedom (DoF)~\citep{Kaneko2015}.
	
	For comparison, we  solve the~\ref{eq:hlsp}'s directly with the ASM solver LexLS (\citep{dimitrov:2015}, https://github.com/jrl-umi3218/lexls), and in sequence~\citep{Kanoun2009} by the quadratic programming (QP) solvers  OSQP~\citep{osqp} based on the alternating direction of multipliers method (ADMM) and the proprietary barrier solver GUROBI~\citep{gurobi}. Unlike our dense solver $\mathcal{N}$IPM-HLSP, the two latter are sparse solvers. 
	Since ADMM based solvers have the tendency to converge with only moderate accuracy, OSQP possesses the functionality to make an educated guess about the active set (`polishing'). We enable this functionality for the second (Sec.~\ref{sec:eval:nonsing}) and third simulation (Sec.~\ref{sec:eval:dynPush}). Otherwise, LexLS (iteration limit: 200), OSQP and GUROBI are used at their standard settings. In case of $\mathcal{N}$IPM-HLSP we limit the number of IPM iterations to 20. The projected normal equations~\eqref{eq:HLSPNeNmethod} are solved by a regularized LDLT decomposition (robust Cholesky decomposition with pivoting~\citep{Bunch1971}, plus diagonal regularization with a small weight of $1\cdot 10^{-6}$). The least-squares form~\eqref{eq:LQNmethod} is solved by the rank revealing QR decomposition.
	
	We depict the overall solver times, the required number of iterations (for OSQP both iterations and number of factorization updates as a dashed line), the time per Newton iteration (solver time divided by number of iterations; not for LexLS and OSQP) and the maximum \ref{eq:kktHLSP} norm of all the levels (not for LexLS since it does not keep an update of the dual corresponding to inactive inequality constraints).
	The simulations are run on an Intel(R) Core(TM) i7-9750H CPU @ 2.60GHz processor with 32 Gb of RAM.
	
	\subsection{Solving equality only LSP's}
	
	\label{eval:eqonly}
	
	\begin{figure}[t!]
		\includegraphics[width=0.6\columnwidth]{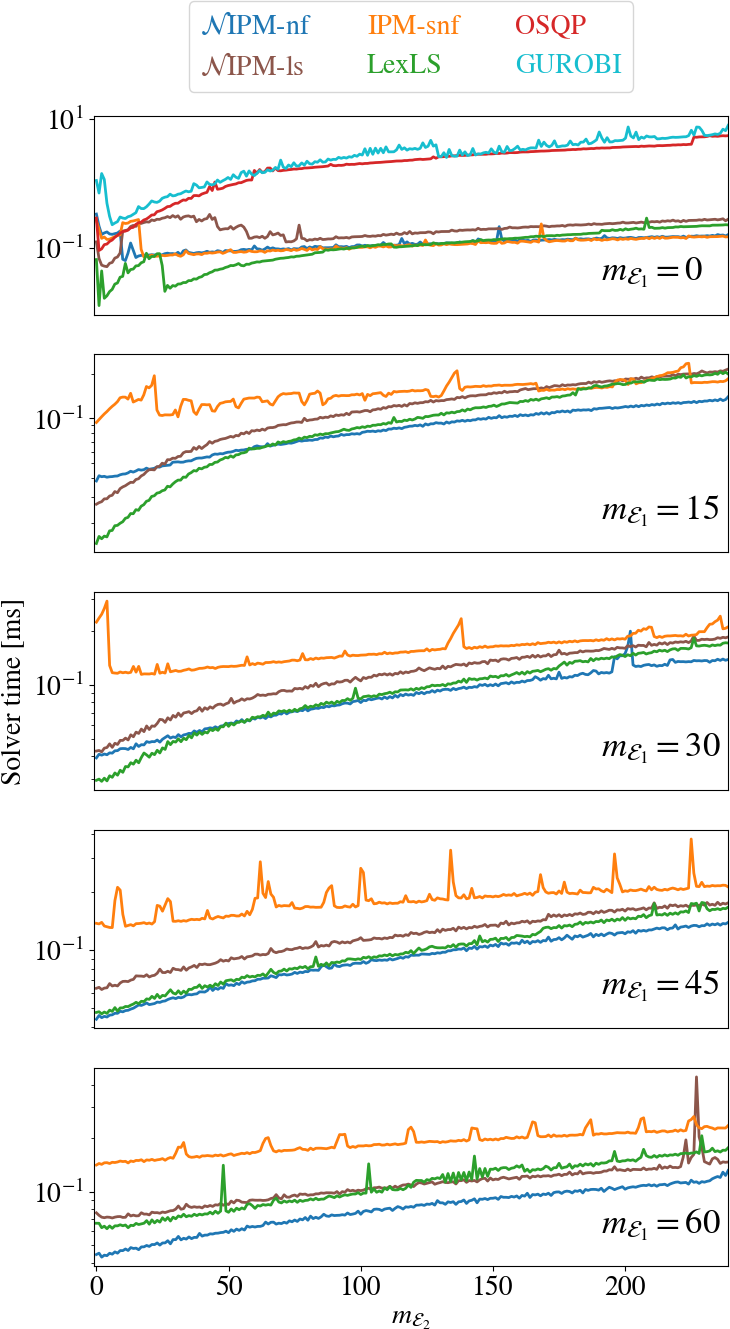}
		\centering
		\caption{Computation times of solving equality only, dense HLSP's ($p=2$). The number of equality constraints on the first level are $m_{\mathbb{E}_1}=0,15,30,45,60$ from top to bottom. The x-axis shows the number of equality constraints on the second level $m_{\mathbb{E}_2}$.}
		\vspace{-10pt}
		\label{fig:benchmark}
	\end{figure}
	
	Figure~\ref{fig:benchmark} shows the computation times for the different solvers when an equality only, dense LSP is solved ($p=2$, $n=60$ variables). The number of constraints of the first level $A_{\mathbb{E}_1}\in\mathbb{R}^{m_{\mathbb{E}_1},n}$ is chosen to $m_{\mathbb{E}_1} = 0,15,30,45,60$. The number of constraints $m_{\mathbb{E}_2}$ on the second level ranges from 0 to 240. 
	
	It can be observed that the problems are solved the fastest for most problems by our solver $\mathcal{N}$IPM-nf.
	The LSP is solved at the same speed by IPM-snf when there are no equality constraints on the first level. In this case both problems are equivalent since IPM-snf only needs to do steps 2 and 6 of Tab.~\ref{tab:NewtonIterOp}. Similarly, $\mathcal{N}$IPM-nf only needs to do steps 2 and 3. 
	For the same problem constellation, both OSQP and GUROBI require much greater computation times due to computational overhead of handling dense problems sparsely.
	Therefore, for $m_{\mathbb{E}_1}>0$ we do not report their results for the sake of better readability. Note that GUROBI includes bound constraints on the problem variables. Therefore any problem becomes an inequality constrained problem which requires more than one solver iteration. This is in contrast to $\mathcal{N}$IPM-nf, $\mathcal{N}$IPM-ls and LexLS which all converge after one iteration for equality only problems.
	
	$\mathcal{N}$IPM-ls and LexLS are slower than $\mathcal{N}$IPM-nf except if both $m_{\mathbb{E}_1}$ and $m_{\mathbb{E}_2}$ are low. At the same time $\mathcal{N}$IPM-ls is slower than LexLS which re-uses the QR decomposition for resolving the primal step to compute the nullspace basis of the active constraints. Additionally, LexLS implements a very efficient in-place block decomposition with minimal memory access overhead.

	As the number of equality constraints $m_{\mathbb{E}_1}$ increases, all nullspace method based solvers $\mathcal{N}$IPM-nf, $\mathcal{N}$IPM-ls and LexLS make up for the computation time of the nullspace basis computation $Z_{\mathbb{E}_1}$ and projection $A_{\mathbb{E}_2}Z_{\mathbb{E}_1}$ as increasingly smaller problems $n_r < n$ need to be solved on the second level. This is in contrast to IPM-snf where an increasingly more expensive second Cholesky decomposition of $A_{\mathbb{E}_1}C_2^{-1}A_{\mathbb{E}_1}$ needs to be computed. Eventually, $n_r=0$ for $m_{\mathbb{E}_1} = n = 60$ which means that $\mathcal{N}$IPM-nf, $\mathcal{N}$IPM-ls and LexLS finish after solving the first level.
	
	This evaluation shows that the projected normal equations~\eqref{eq:HLSPNeNmethod} are computationally equivalent or superior with respect to the Schur normal form~\eqref{eq:qpclassicalNeq} in the limit case of equality constraints only. Since in the presence of inequalities the computational burden of the nullspace method for~\ref{eq:HLSPNeNmethod} will be further offset over the solver iterations, we do not further consider~\ref{eq:qpclassicalNeq} throughout the remainder of the evaluation.

	\subsection{Kinematic robot control}
	\label{sec:eval:clasproj}
	
	\begin{table}[htp!]
		\centering
		\caption{\ref{eq:NL-HLSP}~A with $p=5$ and $n=38$ for the humanoid robot HRP-2.}
		\begin{tabular}{@{} cc @{}}  
			\toprule
			$l$ & Constraint $f_l(u)\geqq v_l$\\
			\midrule
			1 & Lower and upper joint angle limits (76 ineq.)\\
			2 & Two feet and left hand position and orientation (18 eq.)\\
			3 & Center of mass, medial and lateral (2 ineq.)\\
			4 & Right hand position (3 eq.)\\
			5 & Joint velocity regularization (38 eq.)\\
			\bottomrule
		\end{tabular}
		\label{tab:hierarchyA}
	\end{table}

	\begin{figure}[htp!]
		\includegraphics[width=0.8\columnwidth]{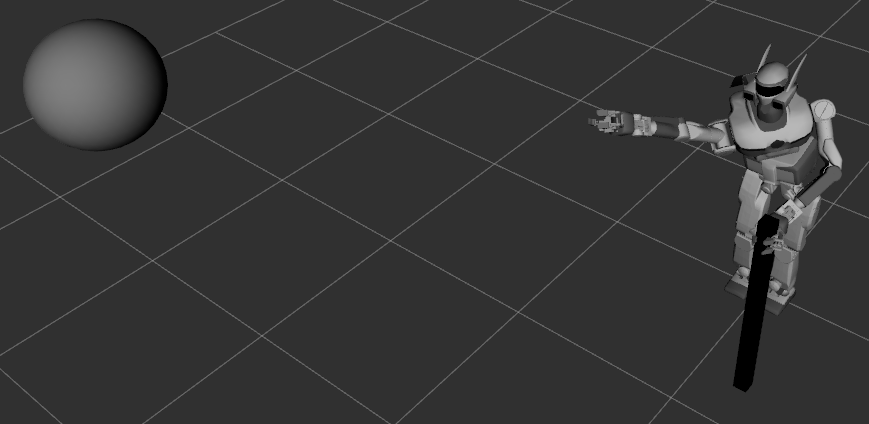}
		\centering
		\caption{HRP-2 stretching into kinematic singularity in order to move the right hand as close as possible to the ball. For the simulation in Sec.~\ref{sec:eval:nonsing}, the robot oscillates in a high frequency and low amplitude manner due to numerical instabilities resulting from unregularized singularities.}
		\label{fig:hrp2_stretch}
	\end{figure}
	
	\begin{figure}[htp!]
		\includegraphics[width=0.8\columnwidth]{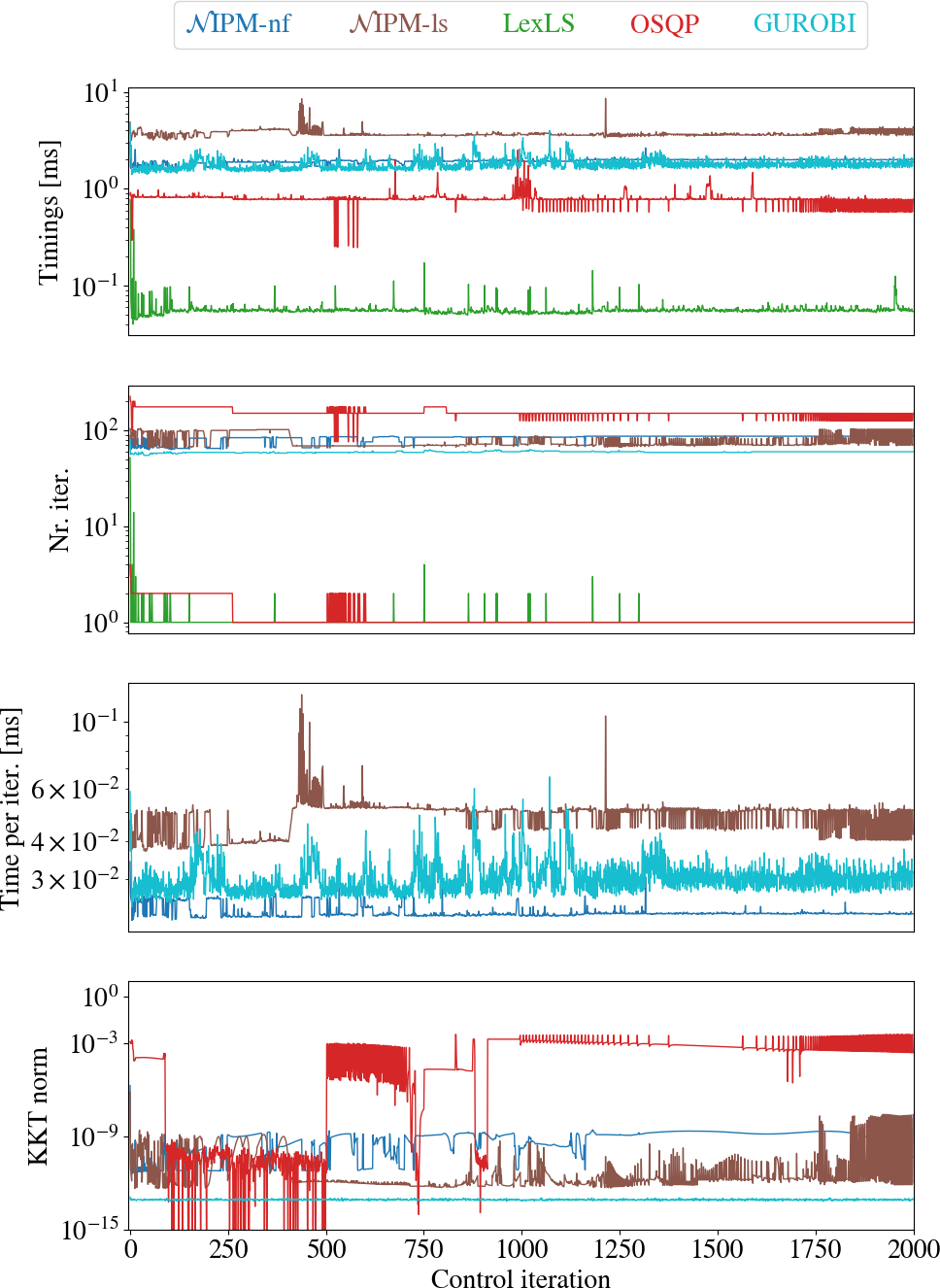}
		\centering
		\caption{\ref{eq:NL-HLSP}~A, Tab.~\ref{tab:hierarchyA}. $\mathcal{N}$IPM-nf exhibits lower computation times per solver iteration than GUROBI.}
		\label{fig:HRP2_2lvl_compClasProjLS}
	\end{figure}
	
	In this example, we demonstrate the computational efficiency of $\mathcal{N}$IPM-HLSP within the S-HLSP solver~\citep{pfeiffer2023} for non-linear real-time robot control problems.
	We design a stretch demonstration with the HRP-2Kai robot according to the control hierarchy~A given in Tab.~\ref{tab:hierarchyA}. The constraints $f$ can represent equality and inequality constraints which is indicated by the symbol $\geqq$. Each constraint's type (inequality, ineq.; equality, eq.) and dimension is given in brackets.
	The first level limits the robot's joint angles. The second level contains kinematic position constraints for both feet to be on the ground and the left hand to be grabbing onto a pole in front of the robot. 
	On the next level the center of mass (CoM) is asked to stay within the 2D projected outer area of the feet on the ground. On the next level the right hand aims to reach an out-of-reach target on the upper right side of the robot. Finally, a regularization task on the joint velocities ensures full rank of the overall hierarchy. This leads to the robot taking the posture depicted in Fig.~\ref{fig:hrp2_stretch}. By virtue of the hierarchical Newton's method~\citep{pfeiffer2023}, the robot control stays numerically stable despite the kinematic singularity of the stretched right arm.
	
	The results are given in Fig.~\ref{fig:HRP2_2lvl_compClasProjLS}. The top graph shows the computation times of the five different solvers in every control time step while the second and third graph show the associated solver iterations and the computation time per solver iteration. The bottom graph shows the sum of the \ref{eq:kktHLSP} norms of each priority level at convergence. 
	
	Among the IPM methods, our solver $\mathcal{N}$IPM-nf has the lowest solver time per solver iteration. The overall solver time is slightly higher than GUROBI since more iterations are required until convergence. All IPM solvers converge to a high degree with the \ref{eq:kktHLSP} norm consistently being lower than $10^{-8}$. This is in contrast to OSQP where in the majority of control iterations the \ref{eq:kktHLSP} norm remains at levels of approximately $10^{-3}$. Note that polishing is disabled as can be seen from the low number of factorization updates (dashed red line in second from top graph). 
	
	The S-HLSP sub-problems are solved fastest by LexLS. After an initial high number of active set changes in the first control iteration (8 active set iterations), the number of active set changes is limited to four, leading to low computation times ($<0.2$~ms) due to warm-starting the active set of each control iteration with the previous one.

	\subsection{Dynamic robot control - A numerically unstable robot movement}
	\label{sec:eval:nonsing}
	
	\begin{table}[htp!]
		\centering
		\caption{\ref{eq:NL-HLSP}~B with $p=6$ and $n=72$ for the humanoid robot HRP-2.}
		\begin{tabular}{@{} cc @{}}  
			\toprule
			$l$ & Constraint $f_l(u)\geqq v_l$\\
			\midrule
			1 & Lower and upper joint angle and velocity limits (152 ineq.)\\
			& Joint torque limits (76 eq.)\\
			& Unilateral contact forces (32 ineq.)\\
			2 & Two feet and left hand position and orientation (18 eq.)\\
			3 & Center of mass, medial and lateral (2 ineq.)\\
			4 & Right hand position (3 eq.)\\
			5 & Joint velocity regularization (38 eq.)\\
			6 & Contact force regularization (36 eq.)\\
			\bottomrule
		\end{tabular}
		\label{tab:hierarchyB}
	\end{table}
	
	\begin{figure}[htp!]
		\includegraphics[width=0.8\columnwidth]{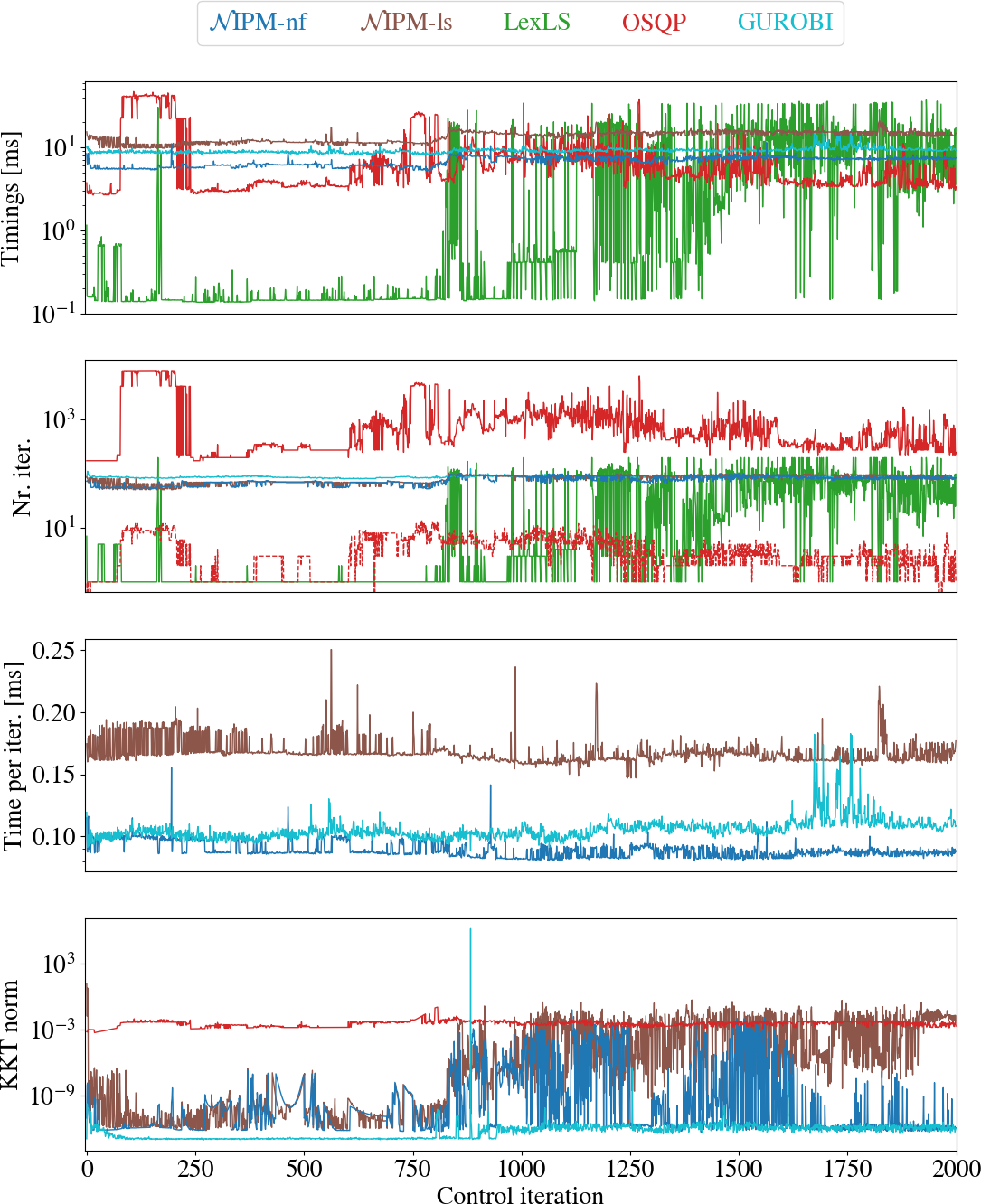}
		\centering
		\caption{\ref{eq:NL-HLSP}~B, Tab.~\ref{tab:hierarchyB}. The right arm starts stretching at control iteration 750 after which the robot behavior becomes numerically unstable.}
		\label{fig:hrp2_dynUnreg}
	\end{figure}

	\begin{figure}[htp!]
		\includegraphics[width=0.3\columnwidth]{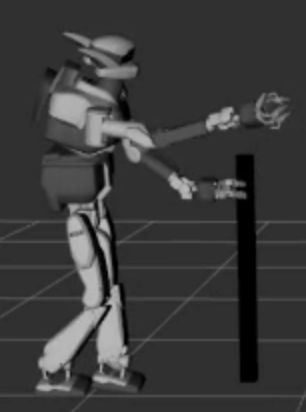}
		\centering
		\caption{\ref{eq:NL-HLSP}~B, Tab.~\ref{tab:hierarchyB}, LexLS. HRP-2 looses contact of the feet with the ground during the stretch task due to sub-optimal active-set identification resulting from constraint singularities. Also, the right arm is not properly stretching towards the goal on the far right top of the robot but swaying around uncontrollably.}
		\label{fig:hrp2ASMfeetloose}
	\end{figure}
	
	\begin{figure}[htp!]
		\includegraphics[width=0.8\columnwidth]{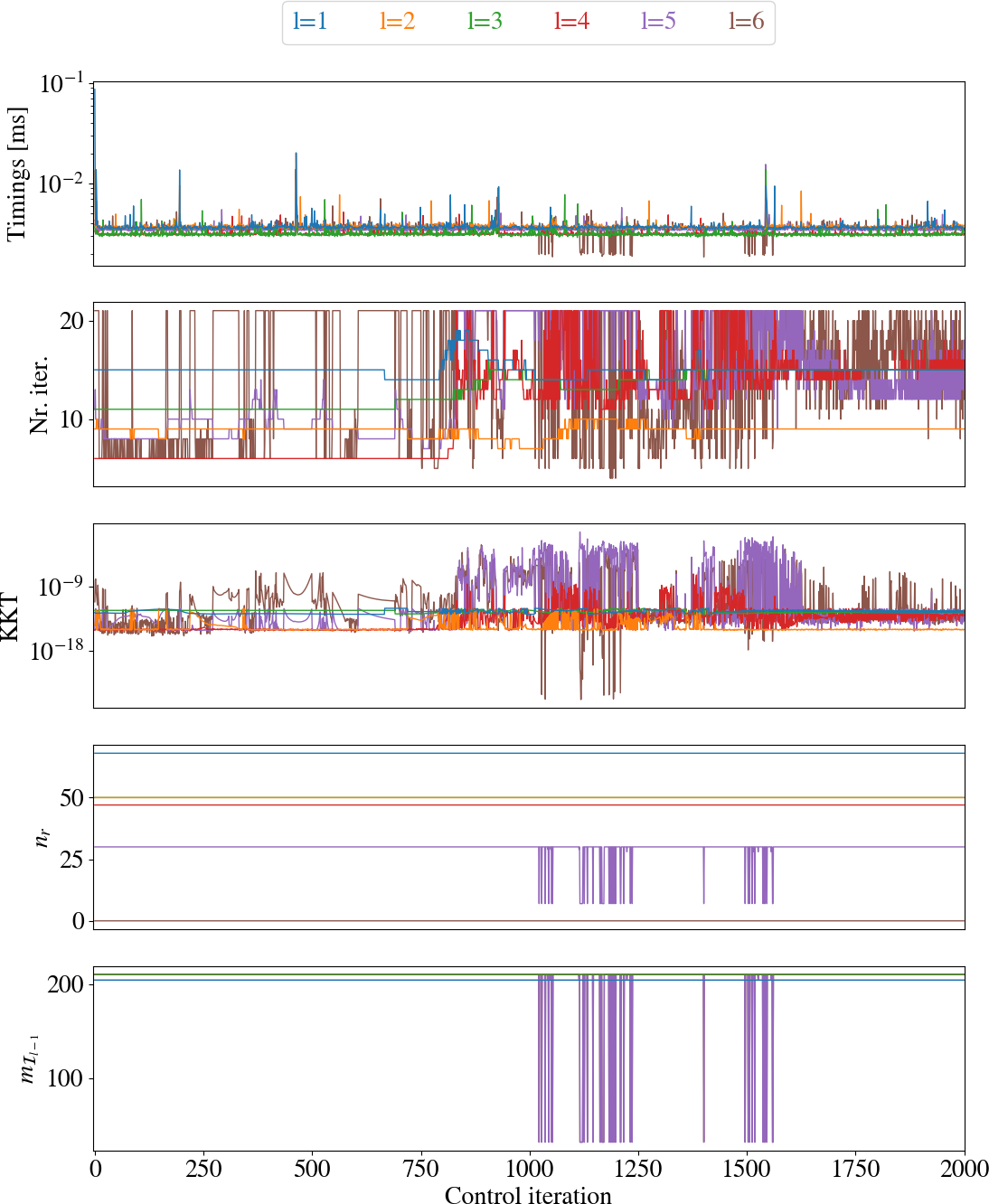}
		\centering
		\caption{\ref{eq:NL-HLSP}~B, Tab.~\ref{tab:hierarchyB}. $\mathcal{N}$IPM-nf's behavior for each priority level. $n_r$ is the number of remaining variables on each respective priority level while $m_{inact}$ corresponds to the number of inactive constraints $m_{\mI_{\cup l-1}}$.}
		\label{fig:hrp2_dynUnreg_NIPMDetails}
	\end{figure}

	The following simulation shows the behavior of the solvers if the problem to solve is scaled up in number of variables and constraints. Additionally, a singular task causes highly dynamic behavior with large changes of the active contact force, joint angle, velocity and torque constraints between the control iterations. The complete control hierarchy B is given in Tab.~\ref{tab:hierarchyB} and the resulting robot posture is given in Fig.~\ref{fig:hrp2_stretch}. We include the equation of motion with 3 contact points on the two feet and the left hand. The joint torques are subsequently substituted by the equation of motion~\citep{Herzog2016} such that we end up with 74 variables consisting of 38 joint velocities and 36 contact forces for the two feet and the left hand.
	The singular task is added to the hierarchy in form of a stretch task trying to reach for a Cartesian position outside of the workspace of the robot. If such a task is not regularized~\citep{Chiaverini1997}, it leads to numerical instabilities in form of high frequency oscillations in presence of a trust region joint velocity constraint~\citep{pfeiffer2023}.
	
	The computation times of the five solvers are given in Fig.~\ref{fig:hrp2_dynUnreg}. $\mathcal{N}$IPM-nf resolves the hierarchy in about 7 ms while the robot is numerically stable. Once the right arm starts stretching at control iteration 750 there is a slight uptick in computation times of about 2 ms due to an increased number of Newton iterations of the levels containing the ill-posed arm stretch task and below ($l \geq 4$), see Fig.~\ref{fig:hrp2_dynUnreg_NIPMDetails}. The fifth level suffers from bad convergence with a maximum \ref{eq:kktHLSP} norm of around $0.06$ ($14.5$ for $\mathcal{N}$IPM-ls). Nonetheless, for both $\mathcal{N}$IPM-nf and $\mathcal{N}$IPM-ls the robot continues to behave as expected with both feet keeping contact to the ground and high frequency oscillations especially on the right arm.
	
	\begin{figure}[htp!]
		\includegraphics[width=0.8\columnwidth]{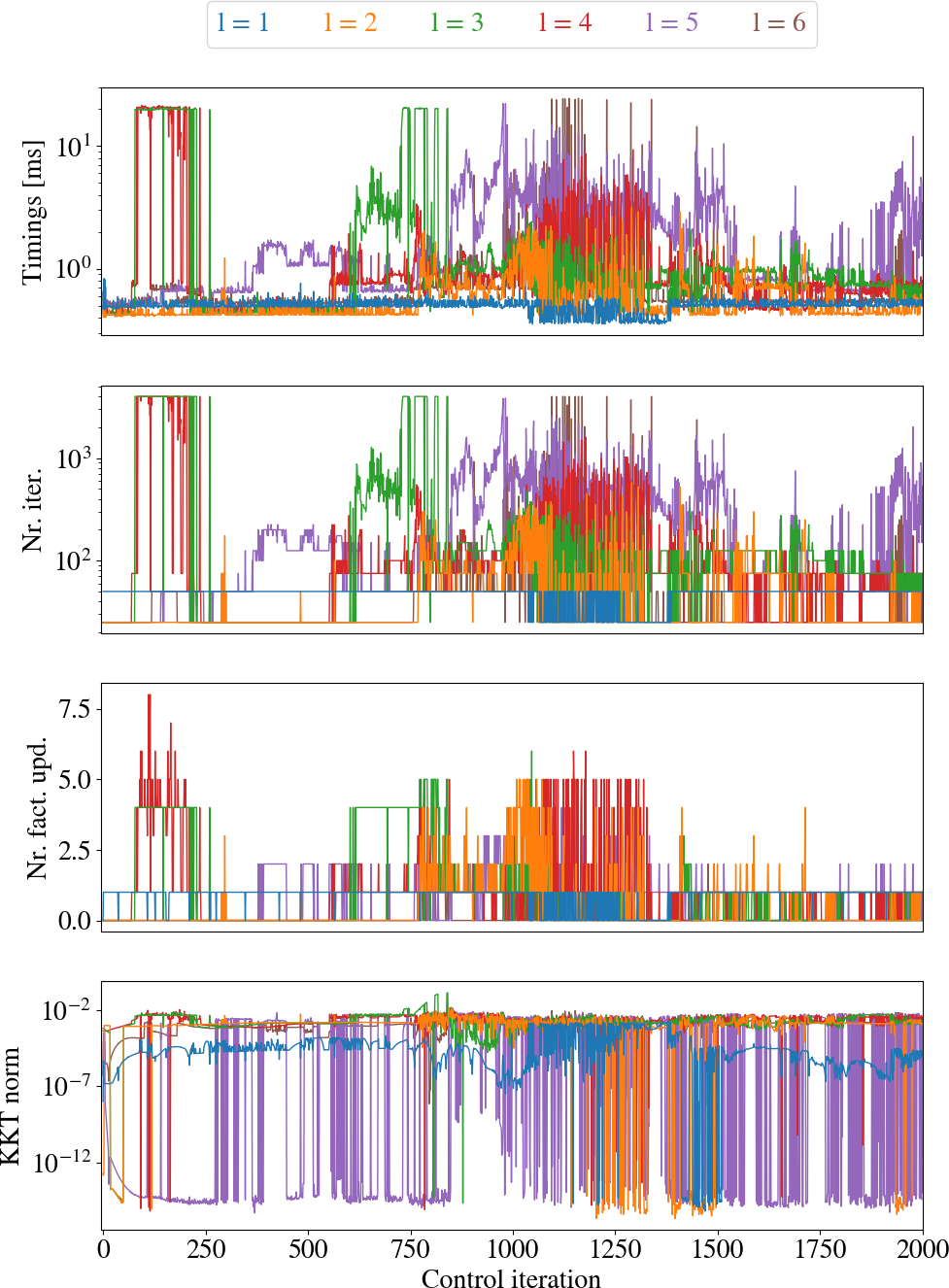}
		\centering
		\caption{\ref{eq:NL-HLSP} B, Tab.~\ref{tab:hierarchyB}. OSQP's behavior for each priority level. `Nr. fact. upd.' (in linear scale) indicates the number of (expensive) factorization updates per priority level.}
		\label{fig:hrp2_dynUnreg_osqpDetails}
	\end{figure}
	
	GUROBI constantly shows a very high convergence rate (except for a single peak of $1.5\times 10^{6}$ due to numerical difficulties). GUROBI efficiently leverages sparsity of the constraint matrices which is not possible for $\mathcal{N}$IPM-HLSP due to dense nullspace projections. Nonetheless, $\mathcal{N}$IPM-nf exhibits approximately the same computation time per solver iteration as GUROBI. This can be explained by the large degree of variable elimination of 32\% from 74 to 50 variables after the first priority level, see second graph from the bottom of Fig.~\ref{fig:hrp2_dynUnreg_NIPMDetails}.
	
	LexLS fails to find a reasonable approximation of the true active set within the given limit of 200 iterations in many control iterations. The high frequency oscillations of the robot due to the kinematic singularities are associated with large changes of the active-set. This leads to the feet losing contact and the right arm swaying uncontrollably, see Fig.~\ref{fig:hrp2ASMfeetloose}. Also, such a large number of active set iterations causes LexLS to violate the real-time constraint with computation times of about 37 ms. 
	
	OSQP converges with a maximum \ref{eq:kktHLSP} norm of about $0.14$ on level 3.
	While the hierarchy is mostly resolved around 2 ms faster than by IPM-ASM, there are many instances where OSQP requires both a high number of iterations and factorization updates. This leads to computation times of up to 42 ms (8175 iterations, 9 factorization updates). We observed that the main challenge of solving~\ref{eq:hlsp} by the ADMM is to maintain feasibility of the subsequent LSP's due to the inherent moderate convergence accuracy of the ADMM. As can be seen from the \ref{eq:kktHLSP} norms in Fig.~\ref{fig:hrp2_dynUnreg_osqpDetails}, there is always a level where the polishing process fails such that the corresponding level converges at about $10^{-3}$. This negatively influences the determination of the optimal infeasible points. While this may be acceptable for LSP's it is problematic for HLSP's since wrongly determined violations $v_{\mA_{\cup l-1}}^*$ render the subsequent priority levels infeasible and eventually lead to solver failure. We explicitly recalculate the violations $v_{\mathbb{E}_l} = A_{\mathbb{E}_l}x - b_{\mathbb{E}_l}$ and $v_{\mathbb{I}_l} = A_{\mathbb{I}_l}x - b_{\mathbb{I}_l}$ at convergence of each priority level which can slightly mitigate this circumstance. Promising progress in terms of accuracy of the ADMM has been documented for example in~\cite{bambade2022}.

	\subsection{Dynamic robot control - Handling a push}
	\label{sec:eval:dynPush}

	\begin{figure}[htp!]
		\includegraphics[width=0.8\columnwidth]{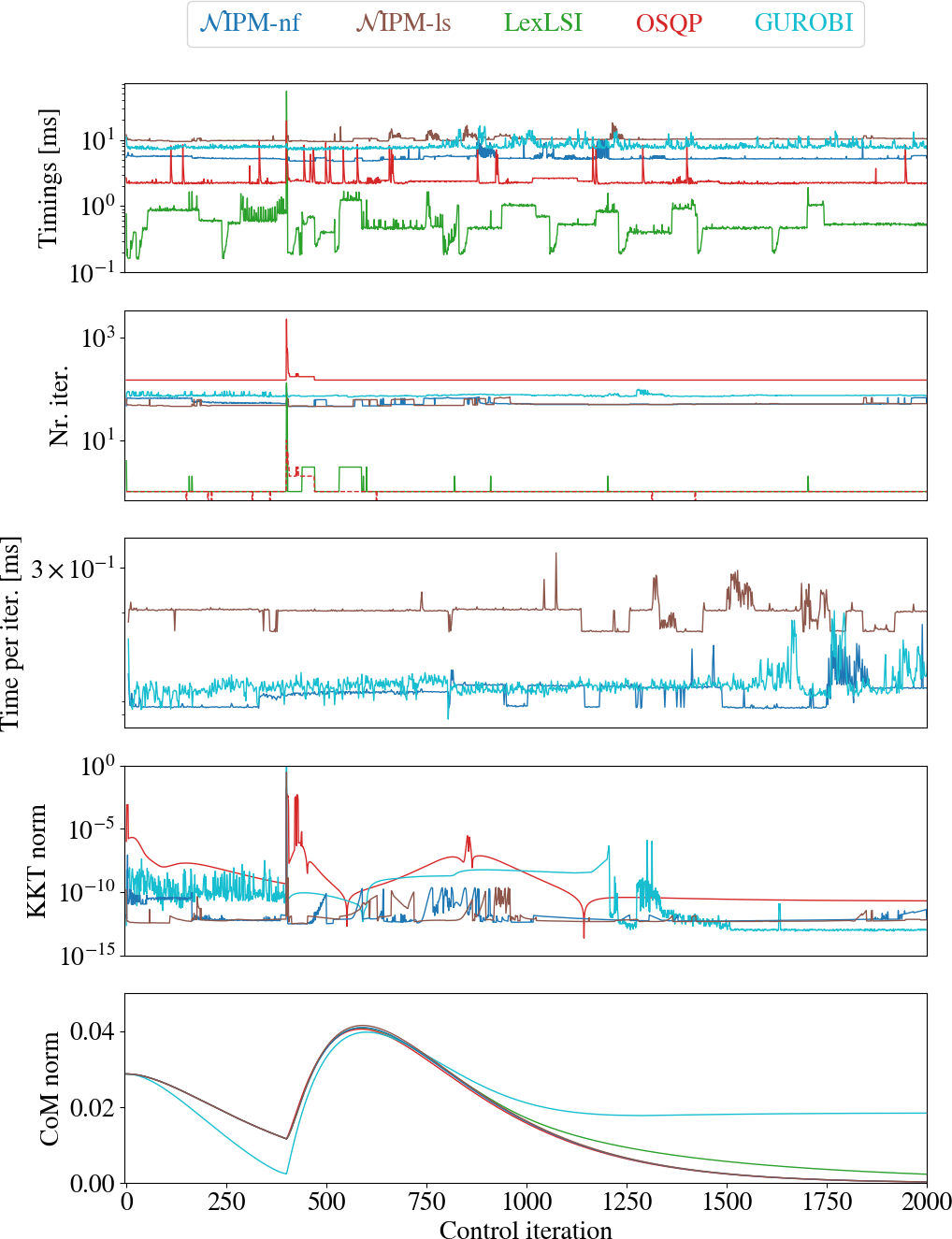}
		\centering
		\caption{\ref{eq:NL-HLSP}~B (Tab.~\ref{tab:hierarchyB}) without right hand task. At control iteration 400 a push from behind of magnitude 400 N is applied to the robot body.}
		\vspace{-10pt}
		\label{fig:hrp2_dynPush}
	\end{figure}
	
	While the previous example in Sec.~\ref{sec:eval:nonsing} is interesting from a numerical point of view, there exist various regularization methods to prevent such numerically unstable behavior~\citep{Chiaverini1997,pfeiffer2023}. This last example therefore looks at a situation which may occur in a real robot situation, namely handling a push. For this we use the control hierarchy B given in Tab.~\ref{tab:hierarchyB} with slight modifications. The level $l=4$ is omitted such that $p=5$. Instead of an inequality constrained CoM, the CoM position is fully controlled by equality constraints in all three spatial directions ($l=3$). 
	
	The robot is asked to lower its CoM. During this process at control iteration 400, a push from behind of magnitude 800 N is applied to the upper body. As can be seen from Fig.~\ref{fig:hrp2_dynPush}, this leads to a singular instance of a large number of 132 active set iterations which takes LexLS approximately 50 ms to resolve. This comes from an unfavorable interplay between the trust region and the joint torque limits caused by possibly ill conditioned equations of motion~\citep{Maciejewski1990}.
	
	This also influences the convergence behavior of OSQP negatively. While OSQP manages successful polishing in most control iterations it does not do so during the push and exposes a higher \ref{eq:kktHLSP} norm in these instances. The high number of iterations and factorization updates require around 19 ms to resolve. 
	
	In contrast, such behavior is observed only to a small degree for the IPM. As can be seen from the second and third graph from the bottom, the number of iterations of $\mathcal{N}$IPM-nf, $\mathcal{N}$IPM-ls and GUROBI are barely influenced. Similarly to the previous example, the \ref{eq:kktHLSP} norm increases sharply in this control instance for all IPM solvers but can be attributed to bad convergence of lower priority levels corresponding to the regularization task. As can be seen from the bottom graph, the CoM trajectory remains physically reasonable as it recovers from the push.
	
	The push leads to the CoM of the robot being shifted to the front and increases the error norm of the CoM task as can be seen from the bottom graph. The difference in CoM behavior seen for GUROBI is due to a regularization term on the joint velocities and contact forces, which amounts to the reformulated cost function $\mini_{x,v_{\mathbb{E}_l},v_{\mathbb{I}_l}} \qquad \frac{1}{2}10^{-5}\Vert x \Vert^2 + \frac{1}{2}\Vert v_{\mathbb{E}_l} \Vert^2 + \frac{1}{2}\Vert v_{\mathbb{I}_l} \Vert^2\qquad$, $l=1,\dots,p$ for the~\ref{eq:hlsp}. Otherwise we observed numerical difficulties.
	
	\section{Conclusion}
	\label{sec:conclusion}
	With this work we have formulated the IPM for HLSP resulting in the solver $\mathcal{N}$IPM-HLSP based on the nullspace method. It requires only a single decomposition of the KKT system per Newton iteration instead of two. This proves computationally equivalent or superior with respect to the IPM in Schur normal form. Our simulations showed that $\mathcal{N}$IPM-HLSP resolves ill-posed problems without significant variations in solver iterations or computation times. In contrast, the ASM tends to require high number of active set iterations in dynamic or numerically unstable control situations due to ill-posed constraint matrices. Our $\mathcal{N}$IPM-HLSP formulation therefore may be preferred to the ASM if solver predictability is regarded more important than very fast computation times but limited by instances of unsuccessful active set searches.
	
	While $\mathcal{N}$IPM-HLSP is reasonably efficient, we see further potential for algorithmic improvements, for example by a sparse solver formulation, heuristically reducing the number of inactive constraints or restricting the number of Newton iterations as seen in~\cite{wangboyd2010}. The last point may require a primal feasible formulation (primal-barrier interior-point method) of our solver. 
	We further believe that our formulation of the nullspace method based IPM for HLSP is not only relevant for instantaneous robotic control but for example in model predictive control (MPC). Special attention requires its block diagonal structure which may be exploited for computational efficiency for example by tailored nullspace bases.
	
	\section{Declarations}
	Part of this work was supported by New York University NSF grants 1925079, 1825993. Part of this work was supported by the Schaeffler Hub for Advanced Research at Nanyang
	Technological University, under the ASTAR IAF-ICP Programme ICP1900093.
	
	\begin{appendices}

		\section{Recursive computation of the Lagrange multipliers associated with active constraints}
		\label{app:recLambdaact}
		\eqref{eq:lambdaact} can be rewritten to
		\begin{align}
			A_{\mA_{\cup l-1}}^T\Delta \lambda_{\mA_{\cup l-1}} = 
			&\BIN
			A_{\mA_{\cup l-1}} \\ A_{\mI_{\cup l-1}} \\ A_{\mathbb{I}_l}  \\ A_{\mathbb{E}_l}  
			\BOUT^T
			\BIN
			-\lambda_{\mA_{\cup l-1}}\\
			D	\\
			E \\
			A_{\mathbb{E}_l}(x+\Delta x) - b_{\mathbb{E}_l}\\
			\BOUT
		\end{align}
		with 
		\begin{align}
			&D \coloneqq -	\lambda_{\mI_{\cup l-1}} -W_{\mI_{\cup l-1}}^{-1}\hspace{-1pt}(\lambda_{\mI_{\cup l-1}} \hspace{-1pt}\odot\hspace{-1pt} ( b_{\mI_{\cup l-1}} \hspace{-3pt}-\hspace{-1pt} A_{\mI_{\cup l-1}}(x+\hspace{-1pt}\Delta x)) + \sigma_{\mI_{\cup l-1}}\mu_{\mI_{\cup l-1}} e)
		\end{align}
		and
		\begin{align}
			E &\coloneqq 
			\left(I + (V_{\mathbb{I}_l} - W_{\mathbb{I}_l})^{-1}{W_{\mathbb{I}_l}}\right)A_{\mathbb{I}_l}\Delta x + A_{\mathbb{I}_l}x - b_{\mathbb{I}_l} - w_{\mathbb{I}_l}  \\
			&+ ({V_{\mathbb{I}_l} - W_{\mathbb{I}_l}})^{-1}({\sigma_{\mathbb{I}_l}\mu_{\mathbb{I}_l} e + w_{\mathbb{I}_l}\odot(A_{\mathbb{I}_l}x - b_{\mathbb{I}_l} - w_{\mathbb{I}_l})})\nonumber
		\end{align}	
		The Lagrange multipliers can be calculated recursively by
		\begin{align}
			&N_{\mA_{\cup j-1}}^TA_{\mathcal{A}_j}^T\Delta \lambda_{\mathcal{A}_j}\hspace{-3pt}=\hspace{-3pt}
			N_{\mA_{\cup j-1}}^T\left(\BIN
			A_{\mI_{\cup l-1}}\\
			A_{\mathbb{I}_l} \\
			A_{\mathbb{E}_l}   
			\BOUT^T
			\hspace{-3pt}
			\BIN
			D\\
			E\\
			A_{\mathbb{E}_l}(x+\Delta x) - b_{\mathbb{E}_l}
			\BOUT \hspace{-3pt}-\hspace{-6pt} \sum_{k=j+1}^{l-1}\hspace{-6pt}A^T_{\mathcal{A}_k}\Delta \lambda_{\mathcal{A}_k}\hspace{-3pt}\right)\label{eq:dualRecursive}
		\end{align}
		with $j=l-1,\dots, 1$.
		For nullspace basis $Z_l$ of the form~\eqref{eq:nsb}, the QR decompositions of $A_{\mathcal{A}_l}N_{\mA_{\cup l-1}}$ can be reused.

		\section{Mehrotra's predictor corrector algorithm for $\mathcal{N}$IPM-HLSP}
		\label{app:predcor}
		
		First, a decomposition of the projected normal equations~\eqref{eq:HLSPNeNmethod} or the least squares form~\eqref{eq:LQNmethod} is computed. Note that this only needs to be done once per Newton iteration. It is then used to first calculate the affine scaling step $\Delta z_{\text{aff}}$, $\Delta x_{\text{aff}}$, $\Delta w_{\mathbb{I}_l}^{\text{aff}}$, $v_{\mathbb{I}_l}^{\text{aff}}$, $\Delta w_{\mI_{\cup l-1}}^{\text{aff}}$ and $\Delta\lambda_{\mI_{\cup l-1}}^{\text{aff}}$ with 
		\begin{align}
			\sigma_{\mI_{\cup l-1}} &= 0 	\quad \text{and} \quad	\sigma_{\mathbb{I}_l} = 0  \label{eq:aff}\\
			F_{\text{aff}} &\coloneqq \lambda_{\mI_{\cup l-1}}
			- W_{\mI_{\cup l-1}}^{-1}\lambda_{\mI_{\cup l-1}} (A_{\mI_{\cup l-1}}x - b_{\mI_{\cup l-1}})\nonumber\\
			G_{\text{aff}} &\coloneqq b_{\mathbb{I}_l} +  w_{\mathbb{I}_l} - A_{\mathbb{I}_l}x -({V_{\mathbb{I}_l} - W_{\mathbb{I}_l}})^{-1}{W_{\mathbb{I}_l}(A_{\mathbb{I}_l}x - b_{\mathbb{I}_l} - w_{\mathbb{I}_l})}\nonumber
		\end{align}
		Line search $\alpha_{\text{aff}}$  is conducted in order to keep the dual feasible with $w_{\mathbb{I}_l} + \alpha_{\text{aff}}\Delta w_{\mathbb{I}_l}^{\text{aff}} \geq 0$, $v_{\mathbb{I}_l} + \alpha_{\text{aff}}\Delta v_{\mathbb{I}_l}^{\text{aff}} \leq 0$, $w_{\mI_{\cup l-1}} + \alpha_{\text{aff}}\Delta w_{\mI_{\cup l-1}}^{\text{aff}} \geq 0$ and $\lambda_{\mI_{\cup l-1}} + \alpha_{\text{aff}}\Delta\lambda_{\mI_{\cup l-1}}^{\text{aff}} \geq 0$.
		
		With this information we calculate the corrector step $\Delta z$, $\Delta x$, $\Delta v_{\mathbb{I}_l}$, $w_{\mathbb{I}_l}$, $\Delta w_{\mI_{\cup l-1}}$ and $\Delta\lambda_{\mI_{\cup l-1}}$ with
		\begin{align}
			\sigma_{\mI_{\cup l-1}} &\coloneqq (\mu_{\mI_{\cup l-1}}^{\text{aff}} / \mu_{\mI_{\cup l-1}})^3 \quad \text{and} \quad\mu_{\mI_{\cup l-1}} \coloneqq \lambda_{\mI_{\cup l-1}}^Tw_{\mI_{\cup l-1}} / m_{\mI_{\cup l-1}}		\label{eq:cntrSigma}\\
			\mu_{\mI_{\cup l-1}}^{\text{aff}} &\coloneqq (\lambda_{\mI_{\cup l-1}} + \alpha_{\text{aff}}\Delta \lambda_{\mI_{\cup l-1}}^{\text{aff}})^T (w_{\mI_{\cup l-1}} + \alpha_{\text{aff}}\Delta w_{\mI_{\cup l-1}}^{\text{aff}}) / m_{\mI_{\cup l-1}}\nonumber\\
			\sigma_{\mathbb{I}_l} &\coloneqq (\mu_{\mathbb{I}_l}^{\text{aff}} / \mu_{\mathbb{I}_l})^3 \quad\text{and}\quad \mu_{\mathbb{I}_l} = -v_{\mathbb{I}_l}^Tw_{\mathbb{I}_l} / m_{\mathbb{I}_l} \nonumber\\
			\mu_{\mathbb{I}_l}^{\text{aff}} &\coloneqq (v_{\mathbb{I}_l} + \alpha_{\text{aff}}\Delta v_{\mathbb{I}_l}^{\text{aff}})^T 
			(w_{\mathbb{I}_l} + \alpha_{\text{aff}}\Delta w_{\mathbb{I}_l}^{\text{aff}})/ m_{\mathbb{I}_l}
			\nonumber
		\end{align}
		For the corrector step we furthermore have
		\\
		\begin{align}
			K_{w_{\mI_{\cup l-1}},l} &\coloneqq v_{\mathbb{I}_l} \odot w_{\mathbb{I}_l}  + \Delta v_{\mathbb{I}_l}^{\text{aff}} \odot \Delta w_{\mathbb{I}_l}^{\text{aff}}  + \sigma_{\mathbb{I}_l}\mu_{\mathbb{I}_l} e \label{eq:cntrK}\\
			K_{w_{\mathbb{I}_l},l}  &\coloneqq \lambda_{\mI_{\cup l-1}}\odot w_{\mI_{\cup l-1}} + \Delta\lambda_{\mI_{\cup l-1}}^{\text{aff}} \odot \Delta w_{\mI_{\cup l-1}}^{\text{aff}} - \sigma_{\mI_{\cup l-1}} \mu_{\mI_{\cup l-1}} e\nonumber
		\end{align}
		such that
		\begin{align}
			D_{\text{cor}} \coloneqq& -	\lambda_{\mI_{\cup l-1}}-W_{\mI_{\cup l-1}}^{-1}(\lambda_{\mI_{\cup l-1}} \odot (b_{\mI_{\cup l-1}}- A_{\mI_{\cup l-1}}(x+\Delta x))\label{eq:cntrD}\\ 
			-& \Delta\lambda_{\mI_{\cup l-1}}^{\text{aff}}\odot \Delta w_{\mI_{\cup l-1}}^{\text{aff}} + \sigma_{\mI_{\cup l-1}} \mu_{\mI_{\cup l-1}} e)\nonumber\\
			F_{\text{cor}} \coloneqq& 	\lambda_{\mI_{\cup l-1}} +W_{\mI_{\cup l-1}}^{-1}(\lambda_{\mI_{\cup l-1}} \odot (b_{\mI_{\cup l-1}} - A_{\mI_{\cup l-1}}x)\nonumber\\ 
			-& \Delta\lambda_{\mI_{\cup l-1}}^{\text{aff}}\odot \Delta w_{\mI_{\cup l-1}}^{\text{aff}} + \sigma_{\mI_{\cup l-1}} \mu_{\mI_{\cup l-1}} e)\nonumber\\
			E_{\text{cor}} \coloneqq& -\left(I + ({V_{1,\mathcal{I}} - W_{1,\mathcal{I}}})^{-1}{W_{\mathbb{I}_l}}\right)A_{\mathbb{I}_l} - A_{\mathbb{I}_l}x + b_{\mathbb{I}_l} + w_{\mathbb{I}_l} \nonumber- ({V_{1,\mathcal{I}} - W_{1,\mathcal{I}}})^{-1}\nonumber\\
			&(\sigma_{\mathbb{I}_l}\mu_{\mathbb{I}_l} e + \Delta v_{\mathbb{I}_l}^{\text{aff}} \odot \Delta w_{\mathbb{I}_l}^{\text{aff}} + w_{\mathbb{I}_l}\odot(A_{\mathbb{I}_l}x - b_{\mathbb{I}_l} - w_{\mathbb{I}_l}) )\nonumber\\
			G_{\text{cor}} \coloneqq& b_{\mathbb{I}_l} + w_{\mathbb{I}_l} - A_{\mathbb{I}_l}x 
			- ({V_{\mathbb{I}_l} - W_{\mathbb{I}_l}})^{-1}\nonumber\\
			&({\sigma_{\mathbb{I}_l}\mu_{\mathbb{I}_l} e + \Delta v_{\mathbb{I}_l}^{\text{aff}} \odot \Delta w_{\mathbb{I}_l}^{\text{aff}} + w_{\mathbb{I}_l}\odot(A_{\mathbb{I}_l}x - b_{\mathbb{I}_l} - w_{\mathbb{I}_l}) })\nonumber
		\end{align}
		
		\section{Algorithm}
		\label{app:alg}
		
		\algnewcommand{\IIf}[1]{\State\algorithmicif\ #1\ \algorithmicthen}
		\algnewcommand{\EndIIf}{\unskip\ \algorithmicend\ \algorithmicif}
		
		\begin{algorithm}[t!]
			\caption{$\mathcal{N}$IPM-HLSP}\label{alg:ipmHLSP}
			\begin{algorithmic}[1]
				\Statex \textbf{Input:} \ref{eq:hlsp}
				\Statex \textbf{Output:} $x$, $\lambda_{\mA_{\cup p}}$
				\State $r = 0$
				\State $\iota = 0$
				\State $N=I$
				\State $\tilde{A} = A_p$
				\For{$l=1:p$}
				\While{$\Vert \tilde{K}_l \Vert_2 > \epsilon$ \& $\iota < \text{maxIter}$}
				\State $\alpha^{\text{aff}},\Delta x^{\text{aff}},\Delta v_{\mathbb{I}_l}^{\text{aff}},\Delta w_{\mathbb{I}_l}^{\text{aff}},\Delta w_{\mI_{\cup l-1}}^{\text{aff}},\Delta\lambda_{\mI_{\cup l-1}}^{\text{aff}}\leftarrow{\tt solve}(\text{`predictor'})$
				\State $\alpha,\Delta x,\Delta v_{\mathbb{I}_l},\Delta w_{\mathbb{I}_l},\Delta w_{\mI_{\cup l-1}},\Delta\lambda_{\mI_{\cup l-1}}$\hspace{-2pt}$\leftarrow$ \hspace{-2pt}${\tt solve}(\text{`corrector'},\alpha^{\text{aff}},\Delta \_^{\text{aff}})$
				\State Make step with $\_ = \_ + \alpha\Delta \_$ for new $x$, $w_{\mathbb{I}_l}$, $v_{\mathbb{I}_l}$, $w_{\mI_{\cup l-1}}$, $\lambda_{\mI_{\cup l-1}}$
				\State $\iota\leftarrow\iota+1$
				\EndWhile\label{euclidendwhile}
				\State $r, N, \tilde{A} \leftarrow {\tt project}(\text{`inactive'}, r, w_{\mI_{\cup l-1}}, \lambda_{\mI_{\cup l-1}}, N, \tilde{A})$
				\IIf{$r \geq n$} \textbf{return} $x$, $\lambda_{\mA_{\cup p}}$\EndIIf 
				\State $r, N, \tilde{A} \leftarrow {\tt project}(\text{`this level'}, r, w_{\mathcal{I}_{l}}, v_{\umI_{l}}, N, \tilde{A})$
				\IIf{$r \geq n$} \textbf{return} $x$, $\lambda_{\mA_{\cup p}}$\EndIIf 
				\EndFor
				\State \textbf{return} $x$, $\lambda_{\mA_{\cup p}}$
			\end{algorithmic}
		\end{algorithm}

		\begin{algorithm}[t!]
			\caption{${\tt solve}$}\label{alg:calc}
			\begin{algorithmic}[1]
				\Statex \textbf{Input:} type, $\alpha^{\text{aff}}$, $\Delta\_^{\text{aff}}$
				\Statex \textbf{Output:} $\alpha$, $\Delta x$, $\Delta v_{\mathbb{I}_l}$, $\Delta w_{\mathbb{I}_l}$, $\Delta w_{\mI_{\cup l-1}}$, $\Delta\lambda_{\mI_{\cup l-1}}$
				\If{type = `predictor'}
				\State Solve~\ref{eq:HLSPNeNmethod}  or~\ref{eq:LQNmethod} for~$\Delta z$ using~\ref{eq:aff} 
				\Else
				\State Solve~\ref{eq:HLSPNeNmethod} or~\ref{eq:LQNmethod} for~$\Delta z$ using~\eqref{eq:cntrSigma},~\eqref{eq:cntrK} and~\eqref{eq:cntrD} with $\alpha^{\text{aff}}$ and $\Delta\_^{\text{aff}}$
				\EndIf
				\State $\Delta x = N_{\mA_{\cup l-1}}\Delta z$
				\State Calculate $\Delta v_{\mathbb{I}_l}$, $\Delta w_{\mathbb{I}_l}$, $\Delta w_{\mI_{\cup l-1}}$, $\Delta\lambda_{\mI_{\cup l-1}}$ with $\Delta x$ 
				\State Line search for $\alpha$ such that~\eqref{eq:linesearch_li} and~\eqref{eq:linesearch_l-1i} are fulfilled
				\State \Return $\alpha$, $\Delta x$, $\Delta v_{\mathbb{I}_l}$, $\Delta w_{\mathbb{I}_l}$, $\Delta w_{\mI_{\cup l-1}}$, $\Delta\lambda_{\mI_{\cup l-1}}$
			\end{algorithmic}
		\end{algorithm}	
		
		\begin{algorithm}[t!]
			\caption{${\tt project}$}\label{alg:project}
			\begin{algorithmic}[1]
				\Statex \textbf{Input:} type, $r$, $w$, $\lambda$, $N$, $\tilde{A}$
				\Statex \textbf{Output:} $r$, $N$, $\tilde{A}$
				\State $\mathcal{A}_l = \{\}$
				\State $\mathbb{I}_l = \{\}$
				\If{type = `inactive'}
				\ForAll{$c \in{\mI_{\cup l-1}}$}
				\If{$w(c) < \xi$ and $\lambda(c) > \xi$}
				\State ${\mathcal{A}_{l^*}} \leftarrow \{c,\mathcal{A}_{l^*}\}$
				\State ${\mI_{\cup l-1}}\leftarrow{\mI_{\cup l-1}} \setminus c$
				\EndIf
				\EndFor
				\State $\hat{r} , Z_{\mA_{l^*}}\leftarrow \mathcal{N}(A_{\mathcal{A}_{l^*}}N_{\mA_{\cup l-1}})$
				\State $N_{\mA_{\cup l^*}} \leftarrow N_{\mA_{\cup l-1}}Z_{\mA_{l^*}}$
				\State $\tilde{A}\leftarrow \tilde{A}Z_{\mA_l^*}$ 
				\Else
				\State ${\mathcal{A}_{l}} \leftarrow \mathbb{E}_l$
				\ForAll{$c \in{\mathcal{I}_{l}}$}
				\If{$w_{\mathcal{I}}(c) < \xi$ and $v_{\mathcal{I}}(c) < -\xi$}
				\State ${\mathcal{A}_{l}} \leftarrow \{c,{\mathcal{A}_{l}}\}$
				\Else
				\State ${\mathcal{I}_{l}}\leftarrow \{c,{\mathcal{I}_{l}}\}$
				\EndIf
				\EndFor		
				\State $\hat{r}, Z_{\mA_l}\leftarrow \mathcal{N}(A_{\mathcal{A}_{l}}N_{\mA_{\cup l^*}})$
				\State $N_{\mA_{\cup l}} \leftarrow N_{\mA_{\cup l^*}}Z_{\mA_{l}}$
				\State  $\tilde{A}\leftarrow \tilde{A}Z_{\mA_l}$
				\EndIf
				\State $r=r + \hat{r}$
				\State \Return $r$, $N_{\mA_{\cup l^*}}$ or $N_{\mA_{\cup l}}$
			\end{algorithmic}
		\end{algorithm}
		
		Here we detail the implementation of $\mathcal{N}$IPM-HLSP which can also be found at \url{https://www.github.com/pfeiffer-kai/NIPM-HLSP}. Algorithm~\ref{alg:ipmHLSP} describes the overall routine. As input it requires a~\ref{eq:hlsp} which contains information about the number of priority levels $p$, the number of variables $n$ and the linear constraints represented by $A$ and $b$. Note that the active and inactive sets $\mA$ and $\mI$ and the optimal slacks $v^*$ are not known yet. Alg.~\ref{alg:calc} summarizes the predictor and corrector step calculation. Alg.~\ref{alg:project} details the active set compositions and projections. The symbol $\_$ is a placeholder for the different variables $x, v_{\mathbb{I}_l}, w_{\mathbb{I}_l},w_{\mI_{\cup l-1}}$ and $\lambda_{\mI_{\cup l-1}}$.
		
	\end{appendices}
	
	\newpage
	
	\bibliography{bib}
	
\end{document}